\documentclass{crmp-l}

\usepackage{amscd,amssymb}
\def\C{{\mathbb C}}
\def\Q{{\mathbb Q}}
\def\R{{\mathbb R}}
\def\Z{{\mathbb Z}}
\def\N{{\mathbb N}}
\def\H{{\mathbb H}}
\def\L{{\mathbb L}}

\def\cA{{\mathcal A}}
\def\cD{{\mathcal D}}
\def\cS{{\mathcal S}}
\def\cC{{\mathcal C}}
\def\cM{{\mathcal M}}
\def\cH{{\mathcal H}}
\def\O{{\mathcal O}}

\def\g{{\mathfrak g}}

\def\k{{\mathfrak k}}
\def\l{{\mathfrak l}}
\def\p{{\mathfrak p}}
\def\u{{\mathfrak u}}

\def\gconj{\overline{g}}

\def\Xconj{\overline{X}}

\def\Dbar{{\overline D}}
\def\Ebar{{\overline E}}
\def\Xbar{{\overline X}}
\def\Ybar{{\overline Y}}
\def\Zbar{{\overline Z}}
\def\Pbar{{\overline P}}

\def\simp{{\boldsymbol{\Delta}}}

\def\-{{-1}}

\def\Hdel{H_{\cD}}
\def\GLn{GL_n(\C)}
\def\gl{{\g\l}}
\def\gln{\gl_n(\C)}

\def\cs{\widehat{c}}
\def\tors{{\mathrm{tors}}}
\def\Hcts{H_{\mathrm{cts}}}
\def\Ccts{C_{\mathrm{cts}}}

\def\dot{{\bullet}}
\def\subdot{_{\dot}}
\def\supdot{^{\dot}}
\def\subdotdot{_{\dot \dot}}

\newcommand\Hom{\operatorname{Hom}}
\newcommand\End{\operatorname{End}}
\newcommand\Spec{\operatorname{Spec}}

\newcommand\cone{\operatorname{cone}}
\newcommand\tr{\operatorname{Tr}}
\newcommand\id{\operatorname{id}}
\newcommand\Ad{\operatorname{Ad}}
\newcommand\image{\operatorname{image}}

\newtheorem{bigtheorem}{Theorem}
\newtheorem{theorem}[equation]{Theorem} 
\newtheorem{lemma}[equation]{Lemma}
\newtheorem{proposition}[equation]{Proposition}
\newtheorem{corollary}[equation]{Corollary}
\newtheorem{statement}[equation]{Statement}
\newtheorem{conjecture}[bigtheorem]{Conjecture}

\newtheorem{bigprop}[bigtheorem]{Proposition}

\theoremstyle{definition}
\newtheorem{definition}[equation]{Definition}
\newtheorem{example}[equation]{Example}
\newtheorem{para}[equation]{}

\theoremstyle{remark}
\newtheorem{remark}[equation]{Remark}

\newtheorem{observation}[equation]{Observation}

\begin{document}

\title[Regulators and Characteristic Classes]{Regulators and
Characteristic Classes of Flat Bundles}

\author{Johan Dupont}

\address{Matematisk Institut\\ Aarhus Universitet\\
DK-80000 \AA rhus C, DENMARK}

\email{dupont@imf.au.dk}

\thanks{First-named author supported in part by grants from
the Aarhus Universitets Forskningsfond, Statens Naturvidenskabelige
Forskningsraad, and the Paul and Gabriella Rosenbaum Foundation.}

\author{Richard Hain}

\address{Department of Mathematics\\ Duke University\\
Durham, NC 27708-0320 USA}

\email{hain@math.duke.edu}

\thanks{Second-named author supported in part by the National Science
Foundation under grants DMS-8601530 and DMS-8901608.}

\author{Steven Zucker}

\address{Department of Mathematics\\ Johns Hopkins University\\
Baltimore, MD 21218 USA}

\email{sz@math.jhu.edu}

\thanks{Third-named author supported in part by the National Science
Foundation under grants DMS-8800355 and DMS-9102233.}
\date{\today}

\begin{abstract}
In this paper, we prove that on any non-singular algebraic variety, the 
characteristic classes of Cheeger-Simons and Beilinson agree whenever
they can be interpreted as elements of the same group (e.g. for flat
bundles). In the universal case, where the base is $BGL(\C)^\delta$,
we show that the universal Cheeger-Simons class is half the Borel
regulator element. We were unable to prove that the universal Beilinson
class and the universal Cheeger-Simons classes agree in this universal
case, but conjecture they do agree.
\end{abstract}


\maketitle

\tableofcontents

\section{Introduction}

For each complex algebraic variety $X$, Deligne and Beilinson \cite{Be}
have defined cohomology groups
$$
H_\cD^k (X, \Z (p)) \quad k, p \in \N
$$
which map to the usual singular cohomology groups:
\begin{equation}\label{1.1}
H_\cD^k (X, \Z (p)) \to H^k (X, \Z(p)).
\end{equation}
(Here, $\Z (p)$ denotes the subgroup $(2\pi i)^p \Z$
of $\C$.)  Beilinson \cite{Be} and Gillet \cite{G} have constructed
Chern classes for algebraic vector bundles $E \to X$
$$
c_p^B (E) \in H_\cD^{2p} (X, \Z (p))
$$
whose images under (\ref{1.1}) are the usual Chern classes. 

In a different
direction, Cheeger and Simons \cite{CS} have defined characteristic
classes
$$
\cs_p (E, \nabla) \in H^{2p-1} (M, \C /\Z (p))
$$
for vector bundles $E \to M$ endowed with a flat connection
$\nabla$ over a smooth manifold.
Flat complex vector bundles over an algebraic manifold $X$ can be made 
algebraic;
we always want to take, as algebraic structure on the bundle, the unique one 
with respect to 
which the connection has regular singularities \cite[p.~98]{De1}.\footnote{If
$E$ is given
a priori as an algebraic bundle on $X$ with flat connection, it is possible 
that the 
two algebraic structures fail to coincide.  See the paragraphs following
\ref{4.2.4}.} 
We can then associate to such bundles the pair of characteristic classes
$$
\cs_p(E, \nabla) \in H^{2p-1} (X, \C /\Z (p)) \text{ and }
c_p^B \in H_\cD^{2p} (X, \Z (p)).
$$

There is a natural homomorphism
\begin{equation}\label{1.2}
H^{2p-1} (X, \C /\Z (p)) \to H_\cD^{2p} (X, \Z(p)).
\end{equation}
Our first result is:

\begin{bigtheorem}
\label{thm1}
If $E \to X$ is a flat algebraic vector bundle with regular singularities over 
the algebraic
manifold $X$, then $c_p^B (E)$ is the image of $\cs_p(E,\nabla)$ under
(\ref{1.2}) for all $p \geq 1$.
\end{bigtheorem}

\noindent Previously known cases of this result are due to Bloch \cite{B}
(flat bundles 
with unitary monodromy over a smooth projective base) and Soul\'e \cite{S} 
(arbitrary 
flat bundles over a smooth projective base).  In \cite{Br}, there is a proof 
of the 
quasi-projective case of Theorem \ref{thm2} below that is similar to ours (in 
\ref{6.2}) and
which invokes our Proposition \ref{6.1.4}.

Since both characteristic classes are functorial, the proof of
Theorem~1 would be rather straightforward if, given $E \to X$ as in the 
theorem,
one could find a morphism of $X$ into some Grassmannian
$f : X \to G (n, \C^m)$ and a connection $\nabla$ on the universal
$n$-plane bundle over $G(n,\C^m)$ that classified {\it both} the bundle
$E$ and its flat connection. However, it is easy to see that this is
impossible in general; the context of flat vector bundles is just too rigid.
For example, if $X$ were compact and $E$ had
nontrivial monodromy, then $f$ could not be constant, implying that
$c_1 (E)$ is non-zero in $H^2 (X, \C)$, which contradicts the
flatness of $E$. As we shall see, there are ways to evade this problem.

We first state a generalization of Theorem~\ref{thm1} that goes beyond the
case of 
flat vector bundles, providing a sort of Chern-Weil theory for DB-cohomolgy 
for algebraic 
manifolds. This applies to algebraic vector bundles with what we call an 
$F^1$-connection --- see \ref{4.2.3}; 
flat 
connections with regular singularities at infinity are examples of 
$F^1$-connections. 
In order to state the result precisely, one has to first write $X=\Xbar -D$,
where $\Xbar$ 
is compact and $D$ is a normal crossings divisor. Then one has to define 
appropriate subgroups
$$
\widehat{H}^{2p}(\Xbar\log D,\Z(p))
$$
of $\widehat{H}^{2p}(X,\Z(p))$ and show that there is a natural map
\begin{equation}
\label{nat_map}
\widehat{H}^{2p}(\Xbar\log D,\Z(p)) \to \Hdel^{2p}(X,\Z(p))
\end{equation}
and that the Cheeger-Simons Chern classes bundles $E$ over $X$ with an
$F^1$-connection lie in this subgroup. This is done in \ref{4.F}.

\begin{bigtheorem}
\label{thm2} If $E \to X$ is an algebraic vector bundle with $F^1$-connection 
over the 
algebraic manifold $X$, then $c_p^B(E)$ is the image of $\cs_p(E,\nabla)$ under
(\ref{nat_map})
for all $p \geq 1$.
\end{bigtheorem}

We first observe that a necessary condition for Theorem~\ref{thm2} is that the
image of 
$\cs_p(E,\nabla)$ in $\Hdel^{2p}(X,\Z(p))$ is the same for all 
$F^1$-connections
$\nabla$ on $E$.  This is verified {\it a priori},---though it appears later
in our 
exposition (as \ref{6.1.4})---and it greatly simplifies the task.

The proof of Theorem \ref{thm2} in the quasi-projective case, given in 
\ref{6.2}, goes 
as follows. For any 
algebraic extension $\Ebar$ of $E$ to a compactification $\Xbar$ of $X$, there 
exists an ample line
bundle $L$ on $\Xbar$ such that both $\Ebar \otimes L$ and $L$ are
generated by global sections, hence, {\it are} pullbacks of universal
bundles by regular maps. Since the conclusion of Theorem~\ref{thm2} is a
tautology for the universal bundles,
we get the result by functoriality for $\Ebar \otimes L$ and $L$; this
implies the same for $\Ebar$ itself, hence for $E$.\footnote{This is,
in essence, the argument in \cite{S}, where $\Xbar = X$; see also our
crucial \ref{6.1.4}.  Note that this argument cannot be used in the
complex analytic setting, not even for compact manifolds.}  The reduction from
algebraic manifolds to 
the quasi-projective case
is fairly standard; see our \ref{40.5}.  

We also formulate and prove the analogous result in the K\"ahler case in 
\ref{10.1}.
In order to effect this, one has to make distinctions among the various
meromorphic 
equivalence classes of compactifications, and then of vector bundle 
extensions.
In the case that $X$ is algbraic, we have been tacitly using the obvious 
choices,
viz.~the algebraic ones (cf.~\ref{conv}).
\medskip

{}From a different point of view, a trick \cite{DK} (see \ref{6.1.5})
can be used to construct a model (depending on $E \to X$ and its
connection) of the universal $n$-plane bundle $E_n \to BGL_n (\C)$ in
the category of {\it simplicial} varieties, a connection on this bundle,
and a morphism
of simplicial varieties $X \to BGL_n(\C)$ that simultaneously classifies
the bundle and the connection. For this reason, it is both natural and
tempting to work in the category of simplicial varieties; Cheeger-Simons
classes, $DB$-cohomology and Beilinson Chern classes all extend directly
to the simplicial setting (cf.\ Appendix~\ref{simp}). We would then want
to view Theorem~\ref{thm2} as a special case of its simplicial analogue.
Unfortunately, there remain difficulties in trying to carry out this
approach; see Section~\ref{discussion}.
An abbreviated account of the above is given in \cite{Z}.

Another reason for working within the simplicial category is that
we can work with the ``universal case'' of a flat bundle.  Denote the
general linear
group of complex $n \times n$  matrices, {\it endowed with the  discrete
topology}, by $\GLn^\delta$.  Its classifying space $BGL_n (\C)^\delta$
classifies flat complex vector bundles of rank $n$.  The Cheeger-Simons
classes of the universal flat bundle
$$
\widehat{c}_p \in H^{2p-1} (BGL_n (\C)^\delta, \C /\Z(p))
$$
are the universal Cheeger-Simons classes.

To generalize a classical theorem of Dirichlet, Borel \cite{Bo} defined
canonical cohomology classes
$$
b_p \in H^{2p-1} (BGL (\C)^\delta, \C/\R(p))
$$
which he used to define ``regulators'' $r_p : K_{2p-1} (\C) \to  \C/\R(p)$
from Quillen's algebraic $K$-theory of
$\C$ into $\C/\R(p)  \cong \Hdel^1(\Spec \C,\R(p))$. Denote by $x_p$ the
class in $H^{2p-1}(\GLn,\Z(p))$ that transgresses to the universal Chern
class $c_p \in H^{2p}(B\GLn,\Z(p))$. For each $n$, the cohomology class
$b_p$ restricts to the image in
$$
H^{2p-1}(\GLn^\delta,\C/\R(p)) \cong H^{2p-1}(B\GLn^\delta,\C/\R(p))
$$
of the element $\xi_p$ of $\Hcts^{2p-1}(\GLn,\C/\R(p))$ that corresponds
to the element $x_p$ of $H^{2p-1}(\GLn,\R(p))$, under the canonical
isomorphism
$$
H^{2p-1}(\GLn,\R(p)) \mathop{\to}\limits^{\cong}
\Hcts^{2p-1}(\GLn,\C/\R(p)).
$$
This construction and the relevant background is reviewed in
Sections \ref{7.1} and \ref{7.2}. Our second main result, proved
in Section~\ref{thm_3} is:

\begin{bigtheorem}
\label{thm4}
The image of $\cs_p$ in $H^{2p-1}(BGL(\C)^\delta, \C/\R(p))$
is $b_p/2$.
\end{bigtheorem}

Since $BGL_n(\C)^\delta$ is a simplicial variety, we also have the
Beilinson Chern classes of the universal flat bundle:
$$
c_p^B \in H_\cD^{2p} (BGL_n (\C)^\delta, \Z (p)) \cong
H^{2p-1}(BGL_n (\C)^\delta, \C /\Z (p)).
$$
The homomorphism $\C/\Z(p) \to \C/\R(p)$ gives rise to a commutative
square:
\begin{equation*}
\begin{CD}
H^{2p}_\cD (BGL (\C)^\delta, \Z (p))  @>{\cong}>>
H^{2p-1} (BGL(\C)^\delta, \C /\Z (p)) \cr
@VVV @VVV \cr
H^{2p}_\cD (BGL (\C)^\delta, \R (p)) @>{\cong}>>
H^{2p-1} (BGL(\C)^\delta, \C/\R(p))
\end{CD}
\end{equation*}

We had hoped to prove that
$$
c_p^B = \cs_p \text{ in } H^{2p-1}(BGL_n(\C)^\delta,\C/\Z(p)).
$$
as a case of the simplicial version of Theorem~\ref{thm2}, which is
stated as \ref{6.1.1}. 
Our approach in trying to prove
that, and some of the difficulties we encountered in our attempt, are 
discussed
in Section~\ref{discussion}. Thus, we are obliged to 
state the preceding as a conjecture:

\begin{conjecture}
\label{conj}
For all $p \geq 1$ and $n \in \N \cup\{\infty \}$,
$$
c_p^B = \cs_p \in H^{2p-1} (BGL_n(\C)^\delta,\C/\Z(p)).
$$
\end{conjecture}

This conjecture, if true, would imply that the universal Beilinson
Chern class for flat bundles is represented by half the
Borel regulator element. Combining Theorem~\ref{thm4} with the
conjecture yields the following conjectural refinement of Beilinson's
result \cite[A5.3]{Be} (see also the article by Rapoport in \cite{RSS})
which asserts that there is a non-zero rational constant $\lambda$ such
that $b_p \equiv \lambda c^B_p$ modulo the decomposable elements
$$
\left[ H^+(\GLn^\delta,\R)\cdot H^+(\GLn^\delta,\R)\right]\otimes \C/\R(p)
\subseteq H\supdot(GL(\C)^\delta,\C/\R(p)).
$$
He used this to prove that his regulators agreed with those of Borel up
to a non-zero rational constant.  Theorem~\ref{thm4} yields as a corollary
the following strengthening of Beilinson's result .

\begin{bigprop} 
\label{thm5}  If Conjecture 4 is true, then the Borel regulator is two times 
the Beilinson
regulator $K_{2p-1} (\C) \to \C/\R(p)$. Consequently, for all number
fields $F$, the Borel regulator and twice the Beilinson regulator
$$
K_{2p-1}(F) \to \Hdel^1(\Spec \, F, \R(p))
\cong \left[\C/\R(p)\right]^{d_p}
$$
are equal.  Here
$$
d_p =
\begin{cases}
r_1 + r_2 & \text{$p$ odd,} \cr
r_2 & \text{$p$ even,}
\end{cases}
$$
where $r_1$ is the number of real embeddings of $F$, and $r_2$ is the
number of conjugate pairs of complex embeddings.
\end{bigprop}

We wish to stress that even if we considered only flat bundles in our
theorems, it would be necessary to consider Cheeger-Simons classes of
bundles with arbitrary connections \cite{CS}, as we make essential use
of the universal $n$-plane bundle which does not admit a flat connection.
By definition, these classes take values in the ring of differential
characters.  We give a detailed account of Cheeger-Simons classes
and differential characters in Section~\ref{cs}, as we could not find
in the literature an exposition suitable for our purposes.  One of the
original approaches to Cheeger-Simons classes uses the construction of a
universal connection by Narasimhan-Ramanan \cite{NR}. In Section~\ref{univ}
we give an alternative and more canonical approach to
universal connections.  In Section~\ref{db_coho}, we review
Deligne-Beilinson cohomology, and prove Theorem~\ref{thm2}. Next, 
Theorem~\ref{thm4} 
is
proved in  Section~\ref{thm_3}.  Finally, we treat the issues and
traps involved in our attempt to prove Conjecture~\ref{conj} in 
Section~\ref{discussion}.

We would like to thank the numerous mathematicians with whom we have had useful
discussions 
related to this paper, especially Pierre Deligne, H\'el\`ene Esnault, Bernard
Shiffman, Vyacheslav Shokurov and 
Christophe 
Soul\'e.  We are particularly grateful to the referee for his or her patience 
over our 
reluctance to recognize the relevance of the distinctions now
made
in \ref{8.9} and \ref{4.F}, even though this was known to the authors.

Hain and Zucker would also like to thank the Max-Planck-Institut f\"{u}r 
Mathematik for 
its hospitality and support during the fall of  1987. Hain would also like to 
thank
the Mathematics Institute at \AA rhus for its support during several visits.

{\numberwithin{equation}{subsection}

\subsection{Conventions}\label{1.3}
To make all Chern classes compatible with those used in algebraic
geometry, we introduce the algebraic geometers' Tate twist. Denote by
$\Z(p)$ the subgroup of $\C$ generated by $(2\pi i)^p$. For each
subgroup $\Lambda$ of $\C$, set
$$
\Lambda (p) = \Lambda \otimes_\Z \Z (p).
$$
The isomorphism $\Z \to \Z (p)$ that takes $1$ to $(2\pi i)^p$
induces a canonical isomorphism
\begin{equation}\label{1.3.1}
H^\dot (X, \Z) \to H^\dot (X, \Z (p)).
\end{equation}
In this paper, the $p^{\rm th}$ Chern class of a $GL_n (\C)$
(equivalently, a complex $n$-plane) bundle over $X$ is the element
of $H^{2p} (X, \Z (p))$ which is the image under (\ref{1.3.1}) of
the usual topological Chern class as defined, for example, in
\cite{MS}.

Let
$$
C_k : {\gl}_n (\C) \to \C \quad k = 0, \ldots, n
$$
be the $GL_n (\C)$-invariant polynomials uniquely determined by
$$
\det (tI - A) = \sum_{k=0}^n C_k (A) t^{n-k}.
$$
Explicitly,
$$
C_k (A) = (-1)^k \tr \wedge^k A.
$$
If $E \to M$ is a complex vector bundle with connection $\nabla$ and
curvature $\Theta$,  then it is well known (see \cite[p.~403]{GH})
that $C_k (\Theta)$ is a closed form that represents $c_k (E)$ in
$H^{2k} (M, \C)$.

All varieties in this paper are defined over the field $\C$ of complex
numbers. We will work in the complex topology unless we explicitly 
say otherwise. The complex of smooth, $\C$-valued forms on a manifold
$X$ will be denoted by $A^\dot(X)$.

\section{Universal Connections}
\label{univ}

\subsection{Classical formulation}\label{2.1}
A construction of universal connections was first given in
\cite{NR}, where the universal bundle $U$ on a Grassmannian $G_\C
(n)$ (or equivalently, its frame bundle) is endowed with a connection
$\nabla^U$, such that any vector bundle with a unitary connection
$(E,\nabla)$ on the manifold $M$ is isomorphic to a pull-back of
$(U,\nabla^U)$ via a classifying mapping $g : M \to G_\C (n)$.
More precisely:

\begin{theorem}\label{2.1.1}{\cite[I\S 4]{NR}}
Fix positive  integers $m$ and $n$. Then there is a positive integer
$\ell$ such that any complex $n$-plane bundle with a unitary connection
$(E, \nabla)$ on an $m$-dimensional manifold $M$ is the pull-back of
the universal one $(U, \nabla^U)$ on the Grassmannian $G(n,\C^\ell)$
of $n$-planes in $\C^\ell$.  Moreover, the  universal connections are
compatible with the inclusions
$G(n,\C^\ell)\hookrightarrow G(n, \C^{\ell+1})$.\endproof
\end{theorem}

\begin{remark}\label{2.1.2}  The integer $\ell$ is given  explicitly in
\cite{NR}. Note that it is necessarily larger than would be needed for
classifying bundles without connection.
Universal connections are constructed for principal bundles with
arbitrary connected structure group in \cite{NR}.
\end{remark}

\subsection{Alternate formulation}\label{2.2}
For our own purposes, it is better to have an alternate formulation
of universal connections. Let $P$ be a principal bundle, with structure
group $G$ (not assumed to  be connected), over the manifold $Y$.
A connection on $P$ is, by definition, a $G$-equivariant lifting of
the tangent bundle $TY$ of $Y$ to $TP$ (the tangent bundle of $P$).
Consider, then, the diagram:
\begin{equation*}
\begin{CD}
TP @>>> TP/G @>p>> TY \cr
@VVV @V{\pi}VV @VVV \cr
P @>>> P/G @= Y
\end{CD}
\end{equation*}
In fact, $\pi$ is a vector bundle projection, and the left-hand
square is cartesian.  Put
\begin{equation}\label{2.2.1}
\widetilde Y = \{ \alpha \in \Hom (TY, TP/G): p \circ \alpha =
\id_{TY}\}
\end{equation}
Let $q : \widetilde Y \to Y$ denote the natural projection. The
following is evident:

\begin{proposition}\label{2.2.2} The connections on $P$ are in
one-to-one correspondence with the cross-sections of $q$.
\end{proposition}

Let $\g$ denote the Lie algebra of $G$. Associated to any connection
is its {\it connection $1$-form}:
$$
\omega \in A^1(P,\Ad(\g)).
$$
In terms of the preceding, $\omega$ can be described as follows.
Let $\widetilde \alpha : TY \to TP$ be the horizontal lift associated
to $\alpha$, and let $\xi \in TP$. Then one has the simple formula:
\begin{equation}\label{2.2.3}
\omega(\xi) =  \xi - \widetilde \alpha (p_*\xi).
\end{equation}

It is convenient to describe $\widetilde Y$ in terms of a local
trivialization of $P$. Thus, we replace $Y$ by a sufficiently small
open subset, {\it which we still call} $Y$.  Then
\begin{equation}\label{2.2.4}
P \,\cong\, Y \times G,
\end{equation}
\begin{equation}\label{2.2.5}
TP/G \,\cong\, TY \times (TG/G) \,\cong \,TY \times {\g},
\end{equation}
so any $\alpha$ in \ref{2.2.1} is determined by an element
$\overline \alpha \in \Hom (TY,{\g})$.  Thus we have:

\begin{proposition}\label{2.2.6}
$\widetilde Y$ is an affine-space bundle over $Y$, with fiber
$\Hom (T_y Y,{\g})$ (for any $y \in Y$).  In particular,
$\widetilde Y$ is a manifold of the same homotopy type as $Y$.
\end{proposition}

Given $f : M \to Y$, the pullback of $P$ is, in terms of \ref{2.2.4},
$$
f^{\ast}P \cong M \times G,
$$
so
\begin{equation}\label{2.2.7}
T(f^\ast P)/G \,\cong\, TM \times {\g}.
\end{equation}
The pullback of the connection is represented by
$\overline \alpha \circ Tf \in  \Hom (TM,{\g})$.

Consider next the diagram
\begin{equation*}
\begin{CD}
@. q^{\ast}TP/G @=  q^{\ast}TP/G @=  q^{\ast}TP/G \cr
@. @AAA @. @VVV \cr
T\widetilde P @>>> T\widetilde P/G @>>> T\widetilde Y @>>>
q^{\ast} TY @>>> TY\cr
@VVV @VVV @VVV @VVV @VVV \cr
\widetilde P @>>> \widetilde P/G @= \widetilde Y
@= \widetilde Y @>q>> Y
\end{CD}
\end{equation*}
where $\widetilde P = q^\ast P$, the pullback
of $P$ along $q$.  By construction, $\widetilde P$
has a tautological connection $\widetilde \nabla$, given by the
pullback of
\begin{equation}\label{2.2.8}
\widehat \beta \in \{ \beta \in \Hom (q^{\ast}TY, q^{\ast}TP/G):
\beta
\text{ projects to } \id_{q^{\ast}TY} \} ,
\end{equation}
where $\widehat \beta$ is defined by: if $\widetilde y \in
q^{-1}(y)$, $\widehat \beta (\widetilde y)$ is
$\widetilde y : T_yY \to (TP/G)_y$.  When $P \cong Y \times G$,
\ref{2.2.8} is determined by
$\overline \beta \in \Hom(q^\ast TY, {\g})$ (recall \ref{2.2.4}
and \ref{2.2.5}).

\begin{corollary}\label{2.2.9}
For any immersion of manifolds $f : M \to Y$ and connection
$\nabla$ on $f^\ast P$, there is a lifting
$\widetilde f : M \to \widetilde Y$ such that
$\widetilde f^\ast \widetilde \nabla = \nabla$.
\end{corollary}

\begin{proof}  By \ref{2.2.2}, the connection $\nabla$
corresponds to a
cross-section $\sigma$ of
\begin{equation*}
\begin{CD}
\widetilde M @= \{\mu \in \Hom (TM, f^\ast TP/G):
f^\ast p \circ \mu = \id_{TM}\}\cr
@VVV \cr M
\end{CD}
\end{equation*}
Since $f$ is, by hypothesis, an immersion, the bundle mapping
$$
Tf : TM \to f^\ast TY
$$
is an injection.  It follows that the natural mapping
\begin{equation*}
\begin{CD}
f^\ast \widetilde Y @= \widetilde Y \times_Y M \subset
\Hom (f^\ast TY, f^\ast TP/G)\cr
@VVrV \cr \widetilde M
\end{CD}
\end{equation*}
is surjective (it is induced by the dual of $Tf$).  Any splitting
of $Tf$ determines a section $s$ of $r$.   Pick one, and take
$\widetilde f = j \circ s \circ \sigma$ (see diagram below).
\begin{equation*}
\begin{CD}
\widetilde M @<r<< f^\ast \widetilde
Y @>j>> \widetilde Y\cr
@VVV @AA{\tilde{f}}A @VVV \cr
M @= M @>>f> Y
\end{CD}
\end{equation*}
\end{proof}

Taking $Y$ to be $G_\C (n)$, and $G = GL_n (\C)$, we
obtain:

\begin{proposition}\label{2.2.10}  The universal bundle
$\widetilde U =  q^\ast U$ on $\widetilde {G}_\C (n)$, with its
tautological connection $\widetilde \nabla$, is universal for
connections on $GL_n(\C)$-bundles.
\end{proposition}

\begin{proof}  Any vector bundle of rank $n$ on a manifold
can be
classified by an immersion into $G_\C (n)$.  Now apply
\ref{2.2.9}.
\end{proof}

{}From the above interpretation of universal connections, one gets,
for free, the following useful fact.

\begin{corollary}\label{2.2.11}
If $\widetilde g_0, \widetilde g_1 : M \to \widetilde G_\C (n)$
are two immersions which classify the bundle with connection
$(E, \nabla)$. Then there is a piecewise-smooth homotopy
$$
\widetilde h : I \times M \to \widetilde G_\C (n).
$$
from $g_0$ to $g_1$ such that $\widetilde h^\ast \widetilde \nabla$
is the ``constant'' connection $d_t \oplus \nabla = p^\ast \nabla$ on
$p^\ast E$ ($p : I \times M \to M$ being the projection).
\end{corollary}

\begin{proof} Let $g_j : M \to G_\C (n)$ be the projection
$q \circ \widetilde g_j$ of $\widetilde g_j$.  As $g_0$ and $g_1$
are necessarily homotopic, pick any homotopy $k$ from $g_0$ to $g_1$;
we may take $k:I \times M \to G_\C (n)$ to be an immersion.  Let
$\widetilde k : I \times M \to \widetilde G_\C (n)$ be a lifting that
classifies $p^\ast (E,\nabla)$.  This is, of course, a homotopy between
its ends,  $\widetilde k_0$ and $\widetilde k_1$.  By construction,
$$
q\circ \widetilde k_j = q\circ \widetilde g_j \quad (j = 0, 1),
$$
so $\widetilde k_j$ and $\widetilde g_j$ can be connected linearly in
the affine-space bundle $\widetilde G_\C (n)$.  These homotopies
also classify the constant connection.  By combining them with
$\widetilde  k$, we obtain $\widetilde h$ as desired.
\end{proof}

\section{The Cheeger-Simons Chern Class}
\label{cs}

In this section we recall and elaborate on the definition and
properties of the Cheeger-Simons invariant, which can be interpreted
as a Chern classes in the ring of differential characters, associated
to a vector bundle with connection on a manifold \cite{CS}.

\subsection{Generalities}\label{3.1}
Let $M$ be a $C^\infty$ manifold, $S_\dot (M)$ the complex
of $C^\infty$ singular chains on $M$ with integer coefficients, and
$Z_\dot (M)$ the subgroup of cycles.  For any abelian group $\Lambda$,
one has
\begin{equation}\label{3.1.1}
S^\dot (M, \Lambda) = \Hom_\Z (S_\dot (M),
\Lambda),
\end{equation}
the corresponding complex of $\Lambda$-valued smooth singular
cochains, whose coboundary operator we shall denote by $\delta$.

Let $A^\dot (M)$ denote the complex of $C^\infty$ differential
forms on $M$ with complex coefficients (one could use real coefficients
here as well). An element $w \in A^k (M)$ determines a cochain
$$
c_w \in S^k (M, \C)
$$
by the formula
\begin{equation}\label{3.1.2}
c_w (\sigma) = \int_{\Delta^k} \sigma^\ast w
\end{equation}
whenever $\sigma : \Delta^k \to M$ is a $C^\infty$ singular simplex.
This embeds $A^\dot (M)$ as a subcomplex of
$S^\dot (M, \C)$.  For any subgroup $\Lambda$ of $\C$,
we then get a morphism of complexes.
\begin{equation}\label{3.1.3}
\iota_\Lambda : A^\dot (M) \to S^\dot (M, \C/\Lambda),
\end{equation}
which is an injection whenever $\Lambda$ is totally disconnected.

\subsection{Differential characters}\label{3.2}
Following \cite{CS}, one makes the following:

\begin{definition}\label{3.2.1} The group of mod~$\Lambda$
{\it differential characters} of degree $k$ on $M$ is the group
$$
\widehat H^k (M, \Lambda) = \{ f, \alpha) \in
\Hom_\Z
(Z_{k-1} (M), \C/\Lambda) \oplus A^k (M): \delta f =
\iota_\Lambda \alpha, \text { and } d \alpha = 0 \}.
$$
\end{definition}

One sees that $\widehat H^\dot (M, \Lambda)$ is a contravariant
functor of $M$; when $\Lambda$ is a ring there is a functorial
product \cite[(1.11)]{CS} that imparts a ring structure to the
differential characters.

\begin{remark}\label{3.2.2} (i) Though \cite{CS} would have us writing
$\widehat H^k (M, \C /\Lambda)$ in \ref{3.2.1}, we have introduced
the above change of notation for the sake of consistency
with that to be used in \S 5. We have also taken the liberty of
shifting by one the degree of a differential character from that of
\cite{CS}, so  that it becomes compatible with other Chern classes.
This also makes the ring structure a graded one.

(ii) If $\Lambda$ is totally disconnected, then $\alpha$ above is
uniquely determined, \ref{3.1.3}, and the second condition of
\ref{3.2.1} is a consequence of the first.
 \end{remark}

It is useful to understand \ref{3.2.1} in terms of conventional
homological algebra.  We have:

\begin{para}\label{3.2.3}
(i)  If $\overline f\in S^{k-1}(M, \C/\Lambda)$ is
any $\Z$-linear extension of $f$ to $S_{k-1} (M)$, then
$\delta \overline f = \iota_\Lambda \alpha$; in particular, if
$\iota (\Lambda) \alpha = 0$, then $\overline f$ is a cocycle.

(ii) If $\overline f_1$ and $\overline f_2$ are two extensions as in
(i), then $\overline f_2 - \overline f_1$ is a $(\C / \Lambda)$-cocycle
vanishing on $Z_{k-1} (M)$ which, by a simple argument using the fact
that $C_{k-1} (M) = Z_{k-1} (M) \oplus F$, where $F$ is free, implies
that the class of $\overline f_2 - \overline f_1$
vanishes in $H^{k-1} (M, \C /\Lambda)$.   We can therefore view
$\widehat H^k (M, \Lambda)$ as a subgroup of
$$
S^{k-1} (M, \C /\Lambda) /\delta S^{k-2} (M, \C/\Lambda)
$$
when $\Lambda$ is totally disconnected.

(iii) The form $\alpha$, representing zero in
$H^k (M, \C /\Lambda)$, has its periods (on $Z_k (M)$) in $\Lambda$.
\end{para}

{}From this, one obtains:

\begin{proposition}\label{3.2.4}
\textup{(i)} There is a canonical and functorial exact sequence
$$
0 \to H^{k-1} (M, \C/\Lambda) \to \widehat H^k (M, \Lambda )
\to A^k_{cl} (M, \Lambda) \to 0,
$$
where
$$
A_{cl}^k (M, \Lambda) = \{ \varphi \in A^k (M) : d\varphi = 0,
\text{ and all periods of $\varphi$  on } Z_k (M) \text{ lie in }
\Lambda \}.
$$

\textup{(ii)} \textup{ (cf.\ \cite[\S 4]{E})}
\begin{align*}
\widehat H^k (M, \Lambda) &\cong H^{k} (\cone
\{ A^{\geq k} (M)
\mathop {\to}\limits^{\iota_\Lambda} S^\dot (M, \C /\Lambda)\} [-1])\cr
&= \H^{k-1} (M, \cone \{ \cA_M^{\geq k} \to \cS_M^\dot(\C /\Lambda) \} );
\end{align*}
here, $\cA_M^\dot$ and $\cS_M^\dot (\C /\Lambda)$ denote
the sheaves of $C^\infty$ forms and singular $(\C /\Lambda)$-cochains
respectively. \endproof
\end{proposition}
\begin{corollary}\label{3.2.5}
If $H^{k-1}(M, \C /\Lambda) = 0$, then 
$$
\widehat H^k (M, \Lambda )\cong A^k_{cl} (M, \Lambda).$$
\end{corollary}

\subsection{Cheeger-Simons classes.}\label{3.3}
Let $(E, \nabla)$ be a $C^\infty$ complex vector bundle with
connection on $M$.  As in Section~\ref{1.3}, the image of the Chern
class
$$
c_p (E) \in H^{2p} (M, \Z (p))
$$
in $H^{2p} (M, \C)$ is represented in de~Rham cohomology by the
polynomial $C_p (\Theta)$ in the curvature $\Theta$ of $\nabla$:
\begin{equation}\label{3.3.1}
C_p (\Theta) \in A^{2p}_{cl} (M, \Z(p)).
\end{equation}
As such, \ref{3.3.1} is functorial for pull-backs of bundles with
connection.  The {\it Cheeger-Simons invariant} of $(E, \nabla)$
is a functorial lifting of $c_p (E, \nabla)$ to
$$
\cs_p (E, \nabla) \in \widehat H^{2p} (M, \Z(p))
$$
in \ref{3.2.4}(i).   In terms of \ref{3.2.4}(ii), it become a
$p^{\rm th}$ Chern class with values in the group of differential
characters with $\Z (p)$ coefficients.

\begin{remark}\label{3.3.2} If $c_p (E, \nabla) = 0$ (e.g., if
$E$ is a flat bundle), then $\cs_p (E, \nabla)$ is an element of
$H^{2p-1} (M,\C /\Z (p))$.
\end{remark}

\subsection{Existence and uniqueness of $\cs_p$}\label{3.4}
The existence and uniqueness of the Cheeger-Simon invariants can be
deduced from the existence of universal connections (\S 2).  By
functoriality, the invariant is completely determined by its value
on the universal connection $(\widetilde U, \widetilde \nabla)$ of
\ref{2.2.10}. On the other hand, because the odd-dimensional
cohomology of a Grassmannian is trivial, there is a unique
lifting $\cs_p (\widetilde U, \widetilde \nabla)$ of its Chern
form to a differential character (see \ref{3.2.4}(i)).  Thus, at
most one Cheeger-Simons invariant can be defined:  if
$\widetilde g : M \to \widetilde G_\C(n)$
is a smooth mapping that classifies $(E, \nabla)$, one {\it must} take
\begin{equation}\label{3.4.1}
\cs_p (E, \nabla) = \widetilde g^\ast \cs (\widetilde
U, \widetilde \nabla).
\end{equation}

One sees that \ref{3.4.1} actually gives a definition of $\cs_p$
by checking that it is, in fact, independent of the choice of $g$.
For this, let $g_0$ and $g_1$ be  two classifying mappings, and
apply \ref{2.2.11} to produce a ``nice'' homotopy $\widetilde h$
between $\widetilde g_0$ and $\widetilde g_1$.  Then the homotopy
formula (see \ref{3.A.13} in Appendix A) yields
\begin{multline}\label{3.4.2}
\widetilde g_1^\ast \cs_p (\widetilde U, \widetilde \nabla)
-
\widetilde g_0^\ast \widetilde c_p (\widetilde U, \widetilde \nabla)
- \delta (B \widetilde h^\ast \cs_p (\widetilde U, \widetilde
\nabla)) \\
 = B (\widetilde h^\ast \delta \cs_p (\widetilde U, \widetilde
\nabla)) = B (\widetilde h^\ast c_p (\widetilde U, \widetilde
\nabla)) = 0,
\end{multline}
for $\widetilde h^\ast c_p (\widetilde U, \widetilde \nabla) \in
A^{2p} (I\times M)$ is pulled back from $M$, so is annihilated by the
fiber integration $B$.

\subsection{An intrinsic construction of $\cs_p$}\label{3.5}
There is a more intrinsic definition of $\cs_p (E, \nabla)$,
one that  uses only functorial constructions on bundles and
connections \cite[\S 4]{CS}. (See also \cite[\S 3]{ch-sim}.)  A truly
intrinsic construction is given in \ref{6.1.5}; that construction uses
simplicial methods.

Given integers $n \geq p \geq 1$, let $V_n^p$ denote the Stiefel
manifold of linearly independent $(n - p+1)$-tuples in $\C^n$.  Then
$V_n^p$ has the homotopy type of its submanifold $U (n)/U (p-1)$ of
unitary $(n-p+1)$-frames. From this, one sees:

\begin{proposition}\label{3.5.1}\textup{(see \cite[(25.7)]{St})}
$$
\pi_i (V^p_n) =
\begin{cases}
0 & \text{for $i < 2p -1$,}\cr
\Z & \text{for $i = 2p-1$}.
\end{cases}
\leqno (i)
$$
\medskip
$(ii)$ $H_{2p-1} (V_n^p, \Z) \cong \Z,$
and is generated by the homology class of
$U (p) / U (p-1) \cong S^{2p-1}$.
\end{proposition}

{}From a complex vector bundle $E$ rank $n$ on $M$, one can form the
associated {\it Stiefel bundle}
$$
\pi : V^p (E) \to M,
$$
with fiber $V_n^p$.  (For instance, $V^1 (E)$ is just the frame
bundle of $E$.)  By construction, $\pi^\ast E$ has $n - p+1$
tautological sections,
i.e., contains a (canonically) trivial bundle of rank $n - p+1$.
Therefore, the Chern classes of $\pi^\ast E$ are those of any
complementary $(p-1)$-plane bundle $W$.   In particular,
\begin{equation}\label{3.5.2}
\pi^\ast c_p (E) =
c_p (\pi^\ast E) = 0 \text{ in } H^{2p}(V^p (E),\Z(p))
\end{equation}
Given any connection $\nabla$ on $E$, the Chern form $\pi^\ast c_p
(E,\nabla)$ is thus exact.  The first goal is to achieve the exactness
in a functorial fashion.

\begin{para}\label{3.5.3}
For any $W$ as above, let $\nabla^W$ denote the unique connection on
$\pi^\ast E$ which satisfies both:
\begin{enumerate}
\item[(i)] $\nabla^W$ and $\pi^\ast \nabla$ agree on $W$,
\item[(ii)] the tautological sections are flat with respect to
$\nabla^W$.
\end{enumerate}
\end{para}

There is a functorial formula expressing the fact that the Chern forms
of two connections represent the same cohomology class
(see \cite[p.~48]{Ch}):
\begin{equation}\label{3.5.4}
c_p (E, \nabla_1) - c_p (E, \nabla_0) = d\eta_p,
\end{equation}
where
\begin{equation}\label{3.5.5}
\eta_p = p \int^1_0 \mu (\omega, \Theta_t, \ldots, \Theta_t)dt;
\end{equation}
here, $\mu$ is the $p$-linear symmetric form with
\begin{equation}\label{3.5.6}
\mu (\Theta, \ldots, \Theta) = C_p (\Theta),
\end{equation}
$\omega = \nabla_1 - \nabla_0$, and $\Theta_t$ is the curvature of
$$
\nabla_t = \nabla_0 + t\omega, ~~t \in [0, 1].
$$

\begin{observation}\label{3.5.7}
Let $\sigma_0, \sigma_1 : M  \to \widetilde M$
be the sections corresponding to $\nabla_0, \nabla_1$ (recall
\ref{2.2.2}); let $h : I \times M \to \widetilde M$ classify
the connection
$$
\overline \nabla : d_t + (1-t) \nabla_0 + t \nabla_1.
$$
Then formula \ref{3.5.4} is just what comes out of the homotopy formula
\ref{3.A.8} when it is applied to the Chern form of the tautological
connection $\widetilde \nabla$ on $\widetilde M$, or equivalently,
when \ref{3.A.7} is applied to the Chern form of $\overline \nabla$.
To see this, one need only observe that the curvature $\Theta$ of $\overline 
\nabla$
is given by
$$
\Theta = \Theta_t + dt \wedge \omega.
$$
\end{observation}

For $\nabla_0 = \nabla^W$ and $\nabla_1 = \pi^\ast \nabla$,
\ref{3.5.5} gives
\begin{equation}\label{3.5.8}
\pi^\ast c_p (E, \nabla) - c_p (\pi^\ast E, \nabla^W)
= d  \eta^W.
\end{equation}
Since
$$
c_p (\pi^\ast E, \nabla^W) = 0 \quad \text { in } A^{2p} (V^p(E)),
$$
we can rewrite \ref{3.5.8} as:
\begin{equation}\label{3.5.9}
\pi^\ast c_p (E, \nabla) = d\eta^W.
\end{equation}
It remains to study the effect of changing $W$.

\begin{proposition}\label{3.5.10} There is a unique natural choice of
$$
\overline\eta (E, \nabla) \in A^{2p-1}(V^p (E))/dA^{2p-2} (V^p(E))
$$
such that  $d\overline \eta (E, \nabla) = \pi^\ast c_p (E, \nabla)$.
\end{proposition}

\begin{proof}  Let $F \cong V_n^p$ denote a fiber of $\pi$. From
\ref{3.5.1} and the Leray spectral sequence for $\pi$, one
deduces the exactness of:
$$
0 \to H^{2p-1} (M, \Z) \mathop {\to}\limits^{\pi^\ast} H^{2p-1}
(V^p (E), \Z) \to H^{2p-1} (F, \Z) \mathop {\to}\limits^\delta H^{2p} (M,
\Z).
$$
The image of one of the generators of $H^{2p-1} (F, \Z) \cong \Z$ under
$\delta$ is,
in fact, $c_p (E)$, which can be seen easily by considering the Serre
spectral sequence of the fibration.  In particular, for the universal
bundle  $U$ on $G_\C (n)$
\begin{equation}\label{3.5.11}
H^{2p-1}(V^p (U), \Z) \cong H^{2p-1}(G_\C (n),\Z) = 0.
\end{equation}
Given $\varphi$, the condition $d \eta = \varphi$ determines $\eta$
up to a {\it closed} form.  From \ref{3.5.11}, a closed
$(2p-1)$-form on $V^p(U)$ is exact.  It follows that
$\overline \eta (\widetilde U, \widetilde \nabla)$ is
uniquely determined.  Therefore, there is at most one natural
construction of $\overline \eta (E, \nabla)$: it must be
$\widetilde g^\ast \overline \eta (\widetilde U, \widetilde \nabla)$
whenever $g : M \to G_\C(n)$ classifies $E$.   The proof that this is
well-defined goes as in \ref{3.4.2}.
\end{proof}

As the preceding is contrary to the spirit with which we began
subsection~\ref{3.5}, we will show directly, after all, that $\eta^W$
(from \ref{3.5.4}) modulo exact forms, independent of $W$. Let $W_0$
and $W_1$ be two such. Consider the following connection on the pullback
of $E$ to $I \times I \times M$:
\begin{equation}\label{3.5.12}
d_s + d_t + (1-t) [(1 -s) \nabla^{W_0} + s \nabla^{W_1}] + t
(\pi^\ast
\nabla);
\end{equation}
it corresponds to a mapping
$$
I \times I \times M \to \widetilde M.
$$
Let
\begin{equation}\label{3.5.13}
B : A^{2p} (I \times I \times M) \to A^{2p-2} (M)
\end{equation}
denote the fiber integration, which is just the iteration of
\ref{3.A.1} in Appendix A; as such, we write
$B = B_2 \circ B_1$.  From \ref{3.A.7}, it follows that
\begin{equation}\label{3.5.14}
d (B \varphi) = (B_2 \varphi_{1, t} - B_2 \varphi_{0, t}) - (B_1
\varphi_{1,s} - B_1 \varphi_{0, s})
\end{equation}
where $\varphi_{0, t}$ is the restriction of $\varphi$ to
$\{ 0 \} \times I \times M$, $\varphi_{1, s}$ to
$I \times \{ 1 \} \times M$, etc.  Let $\varphi$ be the $p$th
Chern form of \ref{3.5.12}.  Then, by \ref{3.5.7}
$$
B_1 \varphi_{1, s} = 0, ~~B_2 \varphi_{1, t} = \eta^{W_1}, \quad
B_2 \varphi_{0, t} = \eta^{W_0}
$$
This gives
$$
d (B\varphi) = \eta^{W_1} - \eta^{W_0} + B_1 \varphi_{0, s}.
$$
The desired conclusion follows from the realization that
$B_1 \varphi_{0, s}$ is the zero form.  The reason for this is that
when $t = 0$, \ref{3.5.12} gives the linear interpolation between
$\nabla^{W_0}$ and $\nabla^{W_1}$, so $\nabla^{W_1} - \nabla^{W_0}$
is an endomorphism of $\pi^\ast E$ that is zero on the tautological
sub-bundle (so is of rank at most $p-1$). The same holds for the
curvature $\Theta_s$.   Since $\mu$ in \ref{3.5.6}
is induced by polarization from trace $\wedge^p$, it follows that
$B_1 \varphi_{0,  s} = 0$.

We can now complete the direct construction of $\cs_p(E,\nabla)$.

\begin{proposition}\label{3.5.15} There is a unique element
$\widehat \eta \in \widehat H^{2p} (M, \Z(p))$
with $\pi^\ast \widehat \eta = [\eta^W,\pi^\ast c_p (E, \nabla) ]$
(and $\delta \widehat \eta =\iota_{\Z (p)} c_p (E, \nabla)$).
\end{proposition}

\begin{proof}  The uniqueness of $\widehat \eta$ is easy,
for it follows from the diagram
\begin{equation*}
\begin{matrix}
\scriptstyle 0 & \scriptstyle \to& \scriptstyle
H^{2p-1} (F, \C /\Z (p)) & \scriptstyle \to& \scriptstyle
\widehat  H^{2p}(F, \Z (p))& \scriptstyle & \scriptstyle &
\scriptstyle & \scriptstyle \cr
& \scriptstyle & \scriptstyle \uparrow& \scriptstyle & \scriptstyle
 \uparrow\cr
\scriptstyle 0& \scriptstyle \to & \scriptstyle
H^{2p-1}(V^p (E),\C /\Z (p))& \scriptstyle \to & \scriptstyle
 \widehat H^{2p}(V^p(E),\Z (p))& \scriptstyle \to & \scriptstyle
A_{cl}^{2p} (V^p (E), \Z (p))& \scriptstyle \to& \scriptstyle 0\cr
& \scriptstyle & \scriptstyle \uparrow& \scriptstyle & \scriptstyle
\uparrow{\scriptstyle \pi^\ast}& \scriptstyle & \scriptstyle \uparrow&
\scriptstyle & \scriptstyle \cr
\scriptstyle 0& \scriptstyle \to& \scriptstyle
H^{2p-1}(M, \C /\Z (p))& \scriptstyle \to& \scriptstyle
\widehat H^{2p} (M,\Z (p))& \scriptstyle \to& \scriptstyle
A^{2p}_{cl}(M, \Z (p))& \scriptstyle \to& \scriptstyle 0\cr
& \scriptstyle & \scriptstyle \uparrow& \scriptstyle & \scriptstyle
& \scriptstyle & \scriptstyle \uparrow& \scriptstyle & \scriptstyle \cr
& \scriptstyle & \scriptstyle 0& \scriptstyle & \scriptstyle &
 \scriptstyle & \scriptstyle 0& \scriptstyle & \scriptstyle
\end{matrix}
\end{equation*}
that $\pi^\ast$ is injective on differential characters.  Likewise,
we argue that $[\eta^W]$ is in the image of $\pi^\ast$ by showing
that the restriction of $\eta^W$ to the fiber $F$ is the zero
differential character, or equivalently,
$$
\int_{S^{2p-1}} \eta^W \in \Z (p).
$$
(Recall that $d \eta^W = \pi^\ast c_p (\nabla, E)$ vanishes on $F$.)
A little bizarrely, we pass over what might have been a direct
calculation on the space $V_n^p$, the partial  frames in the trivial
$n$-bundle over a point, and check it in the universal situation.
There, $S^{2p-1}$ bounds in $V^p (U)$, so we write
$S^{2p-1} = \partial \Gamma$.  Then $Z = \pi_\ast \Gamma$ is a
cycle on the Grassmannian, with
$$
\int_{S^{2p-1}} \eta^W = \int_\Gamma \pi^\ast c_p (\nabla^U) =
\int_Z c_p (\nabla^U) \in \Z (p)
$$
This completes the proof of \ref{3.5.15}.
\end{proof}

With the preceding accomplished, we now make the:

\begin{definition}\label{3.5.16} The differential character
$\widehat \eta$ above is the {\it Cheeger-Simons Chern class}
of $(E, \nabla)$.  We write $\widehat \eta = \cs_p (E, \nabla)$.
\end{definition}

\begin{remark}\label{3.5.17} By \ref{3.2.3}(ii), the Cheeger-Simons
class $\cs_p$ is the choice of an element
$$
y \in S^{2p-1}(M,\C/\Z (p))/\delta S^{2p-2}(M,\C/\Z (p))
$$
with $\delta y = \iota_{\Z (p)} c_p (E, \nabla)$.  In terms of
the cone on $\iota_{\Z (p)}$ (recall \ref{3.2.4}(ii)), it is
represented by
$$
(c_p (E, \nabla), -y) \in A^{2p} (M) \oplus
S^{2p-1} (M, \C /\Z(p))/\delta S^{2p-2} (M, \C /\Z (p)),
$$
which is a cocycle by virtue of the calculation
$$
D (c_p (E, \nabla), -y) =
(-d c_p (E, \nabla), -\delta y + \iota_{\Z (p)} c_p (E, \nabla))
 = (0, 0).
$$
\end{remark}

\section{Deligne-Beilinson Cohomology}
\label{db_coho}

In this section we recall the definition of Deligne-Beilinson cohomology,
and establish the precise relation with differential characters.
References on Deligne-Beilinson cohomology include \cite{Be} and \cite{EV}.

\subsection{Deligne-Beilinson (DB-) cohomology}\label{4.1}
Suppose that $Y$ is a smooth algebraic variety.  By
resolution of singularities \cite{Hi1}, we can find a smooth completion $\Ybar$
of $Y$ such that $\Ybar - Y$ is a normal crossings
divisor $D$.  Denote by $I_p$ the composite
$$
F^p A^\dot (\Ybar \log D) \to
A^\dot Y \mathop {\to}\limits^I S^\dot (Y),
$$
where $F$ denotes the Hodge filtration of the de~Rham complex.

\begin{definition}\label{4.1.1} Suppose that $\Lambda$ is a subring
of $\mathbb R$ and that $p \in \N$.  The {\it Deligne-Beilinson} (or
DB-) {\it cohomology} $H^\dot_\cD (Y, \Lambda (p))$ of
the variety $Y$ is the cohomology of the complex
$$
D^\dot (Y, \Ybar; \Lambda (p)) =
\cone \{ F^p A^\dot ( \Ybar  \log D )
\mathop {\to}\limits^{I_p} S^\dot (Y, \C /\Lambda
(p))\} [-1].
$$
\end{definition}

As the notation suggests, the DB-cohomology of a variety $Y$ is
independent of
the compactification $\Ybar$ chosen.  This follows from
standard arguments (cf.\ \cite[(8.3.2)]{De2}).
As DB-cohomology  is constructed from a cone, we have:

\begin{proposition}\label{4.1.2} For each smooth algebraic variety $Y$, there 
is natural long exact sequence
\begin{multline}\label{les}
\ldots \to H^{k-1}(Y , \C /\Lambda (p)) \to
H_\cD^k (Y, \Lambda (p)) \cr \to F^p H^k( Y ,\C)
\to H^k ( Y, \C / \Lambda (p)) \to \ldots.
\end{multline}
\endproof
\end{proposition}

\begin{corollary}\label{4.1.3} For each variety
\begin{enumerate}
\item[(i)] For each variety $Y$, there is a natural homomorphism
$$
H^{k-1}(Y, \C /\Lambda (p)) \to H_\cD^k (Y,\Lambda (p)).
$$
\item[(ii)] If $H^{2p-1}(Y, \C /\Lambda (p)) = 0$, there is a natural
isomorphism
$$
H_\cD^{2p} (Y,\Lambda (p))\cong H^{p,p}(Y)\cap\, \image\{H^{2p}(Y,\Lambda (p))
\to H^{2p}(Y,\C)\}.
$$
\end{enumerate}
\endproof
\end{corollary}

There are products defined in Deligne cohomology:
\begin{equation}\label{prod}
H_\cD^k (Y,\Lambda (p)) \times H_\cD^l (Y,\Lambda (q))\to H_\cD^{k+l} (Y,
\Lambda (p+q))
\end{equation} 
(see \cite[1.5.1]{Be}). In particular, $\bigoplus _p H_\cD^{2p} (Y,\Lambda 
(p))$ is a ring.

Beilinson \cite{Be} and Gillet \cite{G} have defined Chern classes
$$
c_p^B (E) \in H_\cD^{2p} (Y, \Z (p))
$$
for vector bundles $E \to Y$ over a variety. The construction in 
\cite[1.7]{Be}, 
is in the manner of \cite{Gd}: since the splitting principle 
holds for 
DB-cohomology, the existence of Chern classes in Deligne cohomology reduces to
the
existence of first Chern classes
$$
c_1 (L) \in H^2_\cD (Y, \Z (1))
$$
for line bundles $L \to Y$.  

\subsection{Differential characters and DB-cohomology}
\label{4.2}
Recall from \ref{3.2} that the
mod~$\Lambda$ differential characters of degree $k$ on a $C^{\infty}$
manifold $M$, $\widehat H^k (M, \Lambda)$, are represented by the
cocycles of degree $k-1$ in
$\cone\{ A^{\geq k}(M)\mathop{\to}\limits^{\iota_\Lambda}
S^\dot(M,\C/\Lambda)\}.$

When $X$ is a {\it complex} manifold, of the form $X = \Xbar - D$,
with $D$ a divisor with normal crossings $\Xbar$, we can
incorporate the Hodge filtration and growth conditions
on $A^{\dot}(X)$ to define subgroups of $\widehat H^k (X, \Lambda)$.
Specifically, take
\begin{equation}\label{4.2.1}
F^p\widehat H^k (\Xbar\log D, \Lambda) = H^{k-1}
(\cone \{ F^p A^{\geq k} (\Xbar\log D)
\mathop {\to}\limits^{\iota_\Lambda} S^\dot (X, \C /\Lambda) \} )
\end{equation}
and
$$
\widehat{H}^k(\Xbar\log D,\Lambda) =
F^0\widehat{H}^k(\Xbar\log D,\Lambda).
$$
It is a subgroup of $\widehat H^k (X, \Lambda)$. When $\Xbar$ is compact,
there is an exact sequence:
\begin{equation}\label{4.2.2}
0 \longrightarrow  H^k(X,\C/\Lambda) \longrightarrow
\widehat H^k (\Xbar\log D, \Lambda)
\longrightarrow F^p A_{cl}^k(\Xbar\log D,\Lambda) \longrightarrow 0
\end{equation}
where $A_{cl}^k(\Xbar\log D,\Lambda)$ is the set of closed elements of
$A^\dot(\Xbar\log D)$ whose periods lie in $\Lambda$.

The reason for introducing (\ref{4.2.1}) comes from its similarity
to (\ref{4.1.2}). The cone in (\ref{4.2.1}) is a subset of the 
cone in (\ref{4.1.1}); for $k = 2p$, one obtains, for $X$ algebraic and for
any $\Lambda$, a functorial exact sequence:
\begin{equation}\label{5.5.6}
F^p A^{2p-1}(\Xbar \log D) \rightarrow F^p\widehat H^{2p}
(\Xbar \log D,\Lambda)
\rightarrow H^{2p}_\cD (X, \Lambda) \rightarrow 0.
\end{equation}


\subsection{Meromorphic equivalence}\label{8.9}
Though it may seem premature, it is helpful to leave the setting of algebraic
varieties and 
algebraic vector bundle for awhile.

Let $X$ be a compactifiable complex manifold, in the sense that $X$
admits a
compactification $\overline X$ that is a compact manifold, and 
for which
$D=\overline X - X$ is an analytic subvariety.  We may then modify 
$\overline X$
by blow-ups with smooth center \cite{Hi2} to make $D$ a divisor with 
normal crossings.  

It becomes necessary to divide the compactifications of $X$ into meromorphic 
equivalence classes.  
\begin{definition}\label{1}
Two compactifications $\overline X_1$ and $\overline X_2$ of $X$ are said to 
be {\it 
meromorphically equivalent} if there exists a compactification $\overline X_3$
of $X$
and morphisms of compactifications, i.e., extensions of the identity map of
$X$, 
$\overline X_3\to\overline X_1$ and $\overline X_3\to\overline X_2$.  
\end{definition}
\noindent The above is easily seen to be an equivalence relation.  

We point out that (the underlying complex manifold of) a smooth algebraic 
variety 
can admit compactifications that are not meromorphically equivalent to the 
algebraic 
ones.  A simple example of this is provided by \cite{CE}:
\begin{example}\label{CE}
Let $C$ be an elliptic curve.  Then $X=\C^{p+q}\times C$ admits non-algebraic
compactifications
$\Xbar_{p,q}$ ($p>0$, $0\le q\le p$) with the following properties:
\begin{enumerate}
\item[(i)] There is a principal $C$-fibration $\pi:\Xbar_{p,q}\to\mathbb 
P^p(\C)
\times\mathbb P^q(\C)$.
\item[(ii)] Every meromorphic function on $\Xbar_{p,q}$ is a pullback from
$\mathbb 
P^p(\C)\times\mathbb P^q(\C)$; the algebraic dimension of $\Xbar_{p,q}$
is less
than its complex dimension, so it is not algebraic.
\item[(iii)] $\Xbar_{p,q}$ is diffeomorphic to the product of spheres
$S^{2p+1}
\times S^{2q+1}$.  It follows that $H^2(\Xbar_{p,q},\C) = 0$, and thus 
$\Xbar_{p,q}$
is not K\"ahler.
\end{enumerate}
\end{example}

Furthermore, it is possible for a complex manifold to admit more than one 
algebraic 
structure; it can have inequivalent algebraic compactifications:
\begin{example}\label{2alg}
Let $C$ again be an elliptic curve.  The {\it universal vector extension} of 
$C$, as a
complex manifold, can be given Hodge-theoretically as follows.  Put $H_\Z = 
H^1(C,\Z)$,
$H = H^1(C,\C)$, and $F =  F^1H^1(C,\C)$.  There is a short exact sequence of
abelian groups:
$$
0\to F\to H/H_\Z\to H/(H_\Z + F)\to 0,
$$
which is isomorphic to
$$
0\to \C\to(\C^*)^2\to C\to 0.
$$
The associated $\mathbb P^1$-bundle:
\begin{equation*}
\begin{CD}
\mathbb P^1 @>>> \Xbar \\
@. @VVV \\
{} @. C
\end{CD}
\end{equation*}

\noindent then gives an algebraic variety $\Xbar$ that is a compactification 
of $X = (\C^*)^2$.  
However, $\Xbar$ and $\mathbb P^1\times\mathbb P^1$ are not meromorphically 
(birationally) 
equivalent, for they have non-isomorphic function fields; only the latter one
is an 
algebraic completion of the {\it algebraic} variety $X$.
\end{example}

We also recall the notion of a meromorphic mapping, from complex geometry.
Let $X$ 
and $Y$ be complex analytic varieties, and $\overline X$, $\overline Y$, 
respective 
partial compactifications.
\begin{definition}\label{2}
A map $f:X\to Y$ is said to be {\it meromorphic} with respect to $\overline X$
and $\overline Y$
(and one writes $f:\overline X \dashrightarrow \overline Y$) if the closure of 
the
graph of $f$
in $\overline X \times \overline Y$ is a subvariety of $\overline X \times 
\overline Y$.
\end{definition}
\noindent Because of the existence of resolutions of singularities for complex 
analytic
varieties (\cite{Hi2} again), two compactifications of $X$ are meromorphically 
equivalent 
if and only if the identity map of $X$ defines a meromorphic map between them.

There is an obvious notion of the {\it extendability} of a holomorphic
vector bundle 
$E$ on $X$ to $\Xbar$. 
We next recall the notion of meromorphic equivalence of bundle extensions; it 
is presented in a
somewhat different, though equivalent, manner in \cite [p.~65ff]{De1}).  
\begin{definition}\label{3}
Let $\Xbar$ be a complex 
manifold, 
$D$ a divisor on $\Xbar$, and put $X=\Xbar - D$.  Let $E$ be a vector 
bundle on
$X$.  Two extensions $\Ebar$ and $\Ebar^\prime$ of $E$ to a vector bundle
on $\Xbar$ 
are said to be {\it meromorphically equivalent} if, whenever $U$ is a nice 
open
neighborhood of a point of $D$ in $\Xbar$, on which both bundles are trivial, 
the 
change-of-frame matrix from one to the other (a priori holomorphic on $U\cap 
X$) is 
meromorphic (i.e., does not have essential singularities) along $U\cap D$.
\end{definition}
\begin{remark}\label{merom}
If $\mathfrak M$ is a meromorphic equivalence class of compactifications of 
$X$, there are
corresponding notions of the $\mathfrak M$-extendability of $E$ and 
$\mathfrak M$-meromorphic
equivalence classes of extensions of $E$.
\end{remark}
One can always regard a vector bundle as a complex manifold, forgetting the
linear structure.
The following is almost a tautology:
\begin{lemma}\label{meroeq} With notation as above,
\begin{enumerate}
\item[(i)] Two bundle extensions $\overline E$ and $\overline E^{\,\prime}$ 
are meromorphically
equivalent if and only if $\overline E$ and $\overline E^{\,\prime}$ are 
equivalent 
partial compactifications of $E$.
\item[(ii)] If $\overline E$ and $\overline E^{\,\prime}$ are meromorphically
equivalent, 
then the meromorphic map $\overline E\dashrightarrow \overline E^{\,\prime}$ 
induces
a meromorphic
map $\mathbb P(\overline E) \dashrightarrow \mathbb P(\overline E^{\,\prime})$.
\end{enumerate}
\end{lemma}

\begin{remark}\label{conv}
When $X$ is an algebraic manifold, there is a canonical equivalence class of 
compactifications
of $X$, namely the smooth algebraic completions.  Algebraic vector bundles 
admit algebraic 
extensions to suitable completions, and this class of extensions is likewise 
canonical.
We understand algebraic extensions of algebraic vector bundles on an algebraic
variety 
if it is not specified otherwise.
\end{remark}
\subsection{$F^1$-connections}
\label{4.F}
Throughout, we let $X\subseteq \Xbar$ be as in \ref{4.2}, and $E$ a vector
bundle on $X$.

\begin{definition}\label{4.2.3}
An {\it $F^1$-connection} on $E$ (relative to $\Xbar$) is a connection for 
which there
exists a vector bundle $\Ebar$ over $\Xbar$ such that
\begin{enumerate}
\item[(i)] the restriction of $\Ebar$ to $X$ is $E$;
\item[(ii)] the local connection forms  (with respect to holomorphic
frames of $\Ebar$), and therefore also the curvature, lie in
$F^1 A^\dot (\Xbar\log D, \End(\Ebar))$.
\end{enumerate}
\end{definition}
\noindent We define a sort of ``universal $F^1$-connection'' in Appendix
\ref{4A}.
Note that 
\ref{4.2.3} coincides with \ref{4.2.w} if and only if $X = \Xbar$.

\begin{remark}\label{4.2.4}
This notion occurs in a well-known fact about Hermitian geometry (see 
\cite[p.73]{GH}): 
If $E$ is an Hermitian vector bundle over the {\it compact} complex manifold 
$X$, there 
exists a unique $F^1$-connection on $E$ with respect to which the Hermitian 
metric on $E$ 
is horizontal.  In particular, an extendable vector bundle always admits an 
$F^1$-connection 
(without singularities).  
\end{remark}

The following is evident:

\begin{proposition}\label{4.2.5}
If ${\nabla}$ is an $F^1$-connection on $E$ (relative to $\Xbar$), then 
the Cheeger-Simons class $\cs_p (E, \nabla)$ is in the subgroup
$F^p\widehat H^{2p} (\Xbar\log D, \Lambda)$ of
$\widehat H^{2p} (X, \Lambda)$.
\end{proposition}

To be able to talk sensibly about $F^1$-connections, one needs to know the 
following:
\begin {proposition}\label{mero}
Let $\Xbar$ be a compact complex manifold, $D$ a divisor with normal 
crossings 
on $\Xbar$, and put $X=\Xbar - D$.  Let $(E,\nabla)$ be a vector bundle
with connection on
$X$.  Suppose that $\Ebar$ and $\Ebar^\prime$ are vector bundles on $\Xbar$
with respect
to which $\nabla$ is an $F^1$-connection.  Then $\Ebar$ and $\Ebar^\prime$
are 
meromorphically equivalent.
\end{proposition}
\begin{proof}
The issue is fairly elementary.  Let $U$ be open in $\Xbar$.  If $\omega$ is
the connection
form with respect to a frame of $\Ebar$ on $U\subset \Xbar$, and $\omega^
\prime$ 
is the connection form with respect to a frame of $\Ebar^\prime$ on $U$, 
suppose that
$\omega$ and $\omega^\prime$ have first order-poles\footnote{Since we will
be reducing 
to curves in $U$, there is actually no need to assume that the poles are 
logarithmic,
for the two notions coincide on a curve.} on $U\cap D$.  Let $A$ be the matrix 
expressing 
the frame for $\Ebar^\prime$ in terms of the one for $\Ebar$ (thus $\omega^
\prime = 
A\omega A^\- + (dA)A^\-$).  We wish to conclude that the entries of $A$ must
be meromorphic
along $U\cap D$.

Like holomorphy, meromorphy can be detected by a curve test, i.e., a function 
of several 
variables is meromorphic if and only if its restriction to sufficiently many
curves is---e.g., 
if it is meromorphic in each variable separately---(see \cite{Sh}; 
cf.~\cite[II,\,4.1.1]{De1}).  
Since connections are functorial, we may assume that $U$ is the unit disc in 
$\C$, and $D=\{0\}$.
In this case, the assertion is proved in \cite[II,\,1.19]{De1}.
\end{proof}

\begin{remark}\label{meron} Note that Prop.~\ref{mero} does not say that there 
exists 
$\Ebar$
for which the connection form is logarithmic; nor does it say that if such
$\Ebar$ 
exists, the connection form has logarithmic poles for all $\Ebar^\prime$ 
meromorphically
equivalent to $\Ebar$ (counterexamples abound).  Also, there is no 
contradiction 
between Remark \ref{4.2.4} and Prop.~\ref{mero}.
\end{remark}

According to \cite[II, (5.2)]{De1}, given {\it any} compactification $\Xbar$ 
of $X$, 
a flat bundle $E$ on $X$ admits a vector bundle extension $\Ebar$ to $\Xbar$ 
with respect to which the flat connection can be 
seen
to be an $F^1$-connection.  In fact, according to \ref{mero} any two such 
$\Ebar$ are 
meromorphically equivalent.  We point out that $\Ebar$ can be determined from 
constructions 
that are local on $\Xbar$ along $D$. Local connection forms are then to be 
computed 
on $\Xbar$ in terms of local frames of $\Ebar$.  In particular, there is no 
analogue of
\ref{mero} for equivalence classes of compactifications.

On the other hand, if one starts with an $F^1$-connection on an {\it algebraic}
vector bundle
on $X$, relative to an algebraic compactification $\Xbar$, it does {\it not} 
follow that 
the meromorphic equivalence class of 
extensions of 
$E$ distinguished by \ref{mero} are the algebraic vector bundles on $\Xbar$.
However, 
the two are known to coincide in an important class of examples, namely the 
flat bundles 
``coming from algebraic geometry.''  By that, one means that the fibers of $E$ 
are 
cohomology groups for a family of algebraic varieties over $X$, and $\nabla$ 
is the 
Gauss-Manin connection.  This result is commonly called {\it the regularity 
theorem 
in algebraic 
geometry}.  (See \cite{K} and its generalizations, e.g., \cite[\S 5]{SZ}.)

The following should shed some light on the issue.
\begin{example}\label{exme}
Let $j:X\hookrightarrow\Xbar$ denote the inclusion.  Given the holomorphic 
vector 
bundle $E$, the specification of $\Ebar$ is equivalent to selecting a 
locally-free 
subsheaf ($\O_\Xbar (\Ebar )$) of $j_*\O_X(E)$ on $\Xbar$.  

Suppose that $\dim\, X = 1$ (so $\Xbar$ is automatically algebraic).  Let 
$\Delta\cong 
U \subset \Xbar$ be a disc, with coordinate $t$, for which the
restriction $j|_U$
is $\Delta^*\hookrightarrow\Delta$.  Also, take $E$ to be a line bundle.  
Restricting 
to $U$, we consider the two inequivalent extensions of $E$, $\Ebar = \O_\Delta$
and $\Ebar^
\prime = \O_\Delta\cdot e^{t^\-}$.  Let $\nabla$ be $\tfrac d{dt}$ 
(connection matrix
$\omega = 0$) and $\nabla^\prime = \nabla + \omega^\prime$, where $\omega^
\prime = 
-t^{-2}dt$.  Clearly, $\nabla$ is $F^1$ with respect to $\Ebar$ (indeed
with respect 
to $t^k\Ebar$ for any $k\in\Z$); whereas $\nabla^\prime$, not $F^1$ with 
respect to 
$\Ebar$, {\it is} $F^1$ with respect to $\Ebar^\prime$.  (We leave it to the 
reader to 
squelch the possible misconception that the preceding contradicts GAGA.)  
Note that
the preceding discussion is global when $X=\C ^*$.
\end{example}

Fixing a meromorphic equivalence class $\mathfrak M$ of compactifications of 
$X$, 
we use (\ref{4.2.1}) to define
\begin{equation}\label{5.5.7}
F^p \widehat H^{2p}(X,\Lambda)_{\mathfrak M,\log} = \varinjlim F^p\widehat 
H^{2p}
(\Xbar \log D,\Lambda),
\end{equation}
where the limit is taken over $\Xbar\in \mathfrak M$;
and $F^p A^{2p-1}(X,\,\Lambda)_{\mathfrak M,\log}$ is defined analogously.

The following variant of \ref{4.2.3} is inevitable.
\begin{definition}\label{5.5.9}
An {\it $F^1$-connection} on $E$ {\it (relative to $\mathfrak M$)} is a 
connection
that is $F^1$ relative to some member of $\mathfrak M$.
\end{definition}
\noindent Then \ref{4.2.5} yields immediately:
\begin{proposition}\label{5.5.10}
If ${\nabla}$ is an $F^1$-connection on $E$ relative to $\mathfrak M$,
then the  Cheeger-Simons class
$\cs_p (E, \nabla)$ lies in the subgroup
$F^p\widehat H^{2p} (X, \Lambda)_{\mathfrak M,\log}$ of
$\widehat H^{2p} (X, \Lambda)$.
\end{proposition}

When $X$ is an algebraic manifold and $\mathfrak M$ is the [class generated by]
algebraic completions,
we drop the subscript ``$\mathfrak M$" in (\ref{5.5.7})--(\ref{5.5.10});
cf.~\ref{conv}.

\subsection{Proof of Theorem~\ref{thm2} (quasi-projective case)}
\label{6.2} We return to algebraic manifolds.  First let $X$ be a smooth
projective variety (a fortiori compact), and $E$ an algebraic vector bundle on 
$X$.
There exists an ample line bundle $L$ on $X$ such that both $L$ and
$E' = E \otimes L$ are generated by global sections, so that both
vector bundles are pullbacks of universal bundles on Grassmannians under 
holomorphic
maps. 
We get a pullback diagram:
$$
\begin{matrix}
E & \rightarrow & p_1^*U_n\otimes p_2^*U_1^{-1}\cr
\downarrow &&\downarrow\cr
X & \mathop {\to}\limits^g & G_\C(n) \times G_\C(1)
\end{matrix}
$$
in which $n$ denotes the rank of $E$, $p_1 \circ g$ classifies
$E'$, and $p_2 \circ g$ classifies $L$.

Also, choose any $F^1$-connection $\nabla^L$ on $L$.  The latter,
together with any given $F^1$-connection $\nabla^E$ on $E$,
induces the tensor product $F^1$-connection:\footnote{
in terms of \ref{5.6.2}, the tensor product construction corresponds to
$
\Hom(T^{1,0}Y,\g) \times \Hom(T^{1,0}Y,\C) \to \Hom(T^{1,0}Y,\g).
$}
$$
\nabla^{E'}(e\otimes\ell)=\nabla^E e \otimes\ell+e\otimes\nabla^L\ell.
$$
Likewise, the tensor product of $\nabla^{E'}$ and the dual of
$\nabla^L$ is just $\nabla^E$ again (under the isomorphism
$E\cong E' \otimes L^{-1}$).

As the image in Deligne cohomology is independent of the $F^1$-connection 
(see \ref{6.1.4}), we may assume that the connections are pullbacks of 
$F^1$-connections
on  $G_\C(n)$ and $G_\C(1)$, which exist by \ref{4.2.4}.  Since $G_\C(n)\times 
G_\C(1)$ 
satisfies the hypothesis of \ref{3.2.5} and \ref{4.1.3},\,(ii), the assertion 
of Theorem~\ref{thm2} 
holds for $p_1^*U_n\otimes p_2^*U_1^{-1}$, so by functoriality for $E$. This
completes 
the proof when $X$ is compact.

We also use functoriality to get the assertion for non-compact $X$.  
Specifically, 
let $\Xbar$ be an algebraic completion of the sort considered in \ref{4.2} for 
which $E$ 
extends to an algebraic vector bundle $\Ebar$ on $\Xbar$.  Then there is a 
commutative diagram:
\begin{equation}
\label{quasi}
\begin{CD}
F^p\widehat H^{2p}(\Xbar,\Z(p)) @>>> H^{2p}_{\cD}(\Xbar,\Z(p))\\
@VVV  @VVV\\
F^p\widehat H^{2p}(\Xbar\log D,\Z(p)) @>>> H^{2p}_{\cD}(X,\Z(p))
\end{CD}
\end{equation}
Let $\nabla^\prime$ be any $F^1$-connection on $\Ebar$.  By functoriality,
$\widehat c_p(\Ebar,\nabla^\prime)$ restricts to 
$\widehat c_p(E,\nabla^\prime)$ and $c^B_p(\Ebar)$ maps to $c^B_p(E)$.
By \ref{6.1.4} again,
$\widehat c_p(E,\nabla^\prime)$ has the same image in $H^{2p}_{\cD}(X,\Z(p))$
as $\widehat c_p(E,\nabla)$,
and we are done.
\subsection{General algebraic manifolds}\label{40.5}
We wish to generalize Theorem~\ref{thm2} to the context of general algebraic
manifolds, i.e., to ones that are not quasi-projective, so need not even be
K\"ahler. This is done by
reducing to the quasi-projective case as we next describe.

\begin{proposition}\label{CC}
Suppose that $f:Y\to X$ is a morphism of algebraic varieties satisfying:
\begin{enumerate}
\item[(i)] $f^*: H^{2p-1} (X,\Z)/{tors}\to H^{2p-1} (Y,\Z)/{tors}$ is 
injective 
with image a direct summand in the sense of mixed Hodge structures over $\Z$;
\item[(ii)] $f^*: H^{2p} (X,\Z)_\tors\to H^{2p} (Y,\Z))_\tors$
is injective.
\end{enumerate}
Then $f^*: H^{2p}_{\cD}(X,\Z(p)) \to H^{2p}_{\cD}(Y,\Z(p))$ 
is injective.
\end{proposition}

\begin{proof} The exact sequence \ref{les} gives a commutative diagram with
exact rows:
\begin{equation}
\label{BB}
\begin{CD}
0 @>>> J^{\prime p} (X) @>>> H^{2p}_{\cD}(X,\Z(p)) @>>> F^pH^{2p}(X,\C)\\
@.  @VVV  @VVV  @VVV\\
0 @>>> J^{\prime p}(Y) @>>> H^{2p}_{\cD}(Y,\Z(p)) @>>> F^pH^{2p}(Y,\C)
\end{CD}
\end{equation}
The following is standard:

\begin{lemma}\label{lemm} Given the commutative diagram with exact rows:
$$\CD
0 @>>> A @>>> B @>>> C \\
@.    @V\alpha VV   @V\beta VV   @V\gamma VV \\
0 @>>> A' @>>> B'@>>> C'
\endCD
$$
a sufficient condition for the injectivity of $\beta$ is that $\alpha$ and 
$\gamma$ be injective.
\end{lemma}

We apply the lemma twice. 
The hypotheses imply that the rightmost vertical arrow in \ref{BB} is an 
injection.  Thus, to reach the conclusion,
it suffices by \ref{lemm} to show that leftmost arrow is injective.  Recall
that $J^{\prime p}$ in the above is an extension by the intermediate Jacobian
$J^p$, having the functorial exact sequence
$$
0 \to J^p(X) \to J^{\prime p}(X) \to H^{2p}(X;\Z(p))_\tors \to 0.
$$
The hypothesis implies the injectivity of $J^p(X)\to J^p(Y)$ and of
$H^{2p}
(X;\Z(p))_\tors \to H^{2p}(Y;\Z(p))_\tors$.  By \ref{lemm}, we are done.
\end{proof}

\begin{corollary}\label{FF}
Under the hypotheses of Proposition \ref{CC}, if the 
conclusion of Theorem 2 holds for $Y$, then it also holds for $X$.
\end{corollary}
\begin{proof}
Consider the diagram that one gets by functoriality:
\begin{equation}
\begin{CD}
F^p\widehat H^{2p}(X) @>>> H^{2p}_{\cD}(X,\Z (p))\\
@VVV  @VVV\\
F^p\widehat H^{2p}(Y) @>>> H^{2p}_{\cD}(Y,\Z (p))
\end{CD}
\end{equation}
Then $\cs_p(f^*E,\nabla)\in F^p\widehat H^{2p}(X)$ maps to
$c_p^B (f^*E) = f^*(c_p^B (E))$. Since the right vertical map is
injective by \ref{FF}, $\cs_p(E,\nabla)$ can only
map to $c_p^B (E)$ in $H^{2p}_{\cD} (X,\Z (p))$.
\end{proof}
We next show that such $Y$ exist:
\begin {proposition}\label{DD} 
\begin{enumerate}
\item[(i)] Given a smooth algebraic variety $X$, there is 
a smooth quasi-projective variety $Y$ with surjective birational map $Y\to X$.
\item[(ii)] Under the conditions of (i), there is a direct sum decomposition
in the derived category of $X$:
$$
Rf_*\Z_Y \cong \Z_X \oplus \ker \tr(f).
$$
\item[(iii)] In particular,
$H^\bullet (Y,\Z) \cong H^\bullet (X,\Z)\oplus \ker f_*$,
and the corresponding mixed Hodge structures are isomorphic over $\Z$.
\end{enumerate}
\end{proposition}

\begin {proof}
The argument for (i) was used already in \cite[II:(3.2)]{De2} for extending
Hodge theory to general algebraic varieties, and we sketch  it here.  By a
theorem of Nagata and the Chow Lemma, $X$ fits into a cartesian diagram
\begin{equation}
\begin{CD}
Y& \hookrightarrow & \overline Y\\
@VfVV  @VV{\overline f}V\\
X& \hookrightarrow& \overline X
\end{CD}\label{AA}
\end{equation}
with $\overline X$ a completion of $X$, $\overline Y$ projective, and 
$\overline f$ birational.  (By resolution of singularities \cite{Hi1}, we may 
assume
that $\overline X$ and $\overline Y$ are smooth, and $\overline X - X$ and
$\overline Y - Y$ are divisors with normal crossings.)

Statement (ii) is a topological assertion that follows from the identity
$\tr(f)\circ f^* = \mathbf 1_X$.  Then (iii) follows immediately.
\end{proof}

Combining Propositions \ref{CC} and \ref{DD}, and keeping in mind our 
convention \ref{conv}, we obtain:

\begin{theorem}\label{2a}
If $X$ is a non-singular complex algebraic variety, $E$ an algebraic vector 
bundle on 
$X$ that extends to an algebraic vector bundle on some completion of $X$, and 
$\nabla$  
an $F^1$-connection on $E$ (relative to said extension of $E$), then
$c_p^B (E)$ 
is the image of $\cs_p(E,\nabla)$ for all $p \geq 1$.
\end{theorem}

\begin{corollary}
The conclusion of Theorem~\ref{thm1} holds for all non-singular algebraic
varieties.
\end{corollary}

\subsection{Chern classes in $\mathfrak  M$-DB-cohomology}\label{10.1}
We wish to define the analogue of $c_p^B(E)$ in the complex analytic setting. 
Fix the  meromorphic equivalence class $\mathfrak M$ of compactifications 
$\overline X$ of $X$. We can use the formula in Def.~\ref{4.1.1} to define
$H^\dot_\cD (X,\Xbar;  \Lambda (p))$.
When $\Xbar$ is K\"ahler, 
this is seen to be independent of $\Xbar\in\mathfrak 
M$, for the same reason it was true in 
\ref{4.1} (see \ref{hodge}), so we rename it:
\begin{definition}\label{MDB}
The $\mathfrak M$-DB-cohomology $H^\dot_\cD (X,\Lambda (p))_{\mathfrak M}$
of $X$ with coefficients in $\Lambda(p)$ is the common 
value of $H^\dot_\cD (X,\Xbar;\Lambda (p))$ for $\Xbar\in\mathfrak M$.
\end{definition}
Then (\ref{5.5.6}) gives rise to the exact sequence:
\begin{equation}\label{5.5.8}
F^p A^{2p-1}(X,\,\Lambda)_{\mathfrak M,\log}\rightarrow F^p \widehat H^{2p}(X,
\Lambda)_{
\mathfrak M,\log}\rightarrow H^{2p}_\cD (X,\Lambda)_{\mathfrak M}\rightarrow 0.
\end{equation}

Next, fix a meromorphic equivalence class $\mathfrak E$ of extensions of $E$
to the members 
of $\mathfrak M$ (note that $\mathfrak E$ subsumes $\mathfrak M$).  We call 
these the $\mathfrak E$-extensions
of the 
vector bundle $E$ on $X$.  We will construct Chern classes 
\begin{equation}\label{CMBD}
c^B_p(E;\mathfrak E)\in H^\dot_\cD (X,\Z (p))_{\mathfrak M}
\end{equation}
in more or less the same way as in the algebraic setting.  Choose $\Xbar\in
\mathfrak M$ to 
which $E$ 
extends to $\Ebar\in\mathfrak E$, and let $j:X\hookrightarrow\Xbar$ denote the
inclusion.

As usual, we consider first the case where $E$ is a line bundle. 
Then there is a short exact sequence of complexes of sheaves on $\overline X$:
\begin{multline}\label{DBS}
0\to \cone\{\Z_{\overline X}(1)\to \mathcal O_{\overline X}\}[-1]
\to \cone\{Rj_*\Z_X(1)
\to \mathcal O_{\overline X}\}[-1]\cr
\to \cone\{\Z_{\overline X}(1)\to Rj_*\Z_X(1)\}\to 0.
\end{multline}

\noindent From this, we extract the following commutative diagram:
\begin{equation*}\label{DBH}
\begin{CD}
\Hdel^2(\overline X,\,\Z(1))  @>>> \Hdel^2(X,\,\Z(1))_{\mathfrak M} @. {} \\
 @A\simeq A d\log A  @AA c^B_1(\mathfrak E) A  @. \\
H^1(\overline X,\mathcal O^*_{\overline X})  @>>>  H^1(X,
\mathcal O^*_X)_{
\mathfrak E}\,\, @. \hookrightarrow H^1(X,\mathcal O^*_X),
\end{CD}
\end{equation*}
\noindent where $H^1(X,\mathcal O^*_X)_{\mathfrak E}$ denotes the 
$\mathfrak E$-extendable line bundles on $X$.
Since any two extensions in $\mathfrak E$ differ by a divisor class supported
on $D$, we see that the Chern class $c_1^B(E;\mathfrak E)\in 
\Hdel^2(
X,\,\Z(1))_{\mathfrak M}$ of an $\mathfrak E$-extendable line 
bundle $E$
is well-defined.

To obtain the same for higher-rank $E$, one invokes the splitting 
principle.   
In the formulation of \cite[p.\,140,\,A1]{Gd}, one must know \ref{split} 
below.
Because we will invoke a little Hodge theory in the argument, we must impose
a condition on $X$.
\begin{definition}\label{qck}
A complex manifold $X$ is said to be {\it quasi-c-K\"ahler} if it admits a
compactification 
$\overline X$ that is a compact K\"ahler manifold, and for which $D=\overline 
X - X$ is an 
analytic subvariety.
\end{definition}
When $X$ is quasi-c-K\"ahler, one can then modify such $\overline X$ by a 
sequence 
of blow-ups with smooth center, a process that preserves the K\"ahler 
condition, to make $D$ a divisor with normal crossings \cite{Hi2}.  

Let
$\psi:\mathbb P(E)\to X$ be the projectivization of $E$, $L$ the
tautological line sub-bundle $\O(1)$ of $\psi^*(E)$, and
$\widetilde{\mathfrak M}$ the meromorphic equivalence 
class of
$\{\mathbb P(\Ebar): \Ebar\text{ an $\mathfrak E$-extension of } E\}$.
\begin{proposition}\label{split}
Let $X$ be quasi-c-K\"ahler, and $\mathfrak M$ a meromorphic equivalence class
of K\"ahler 
compactifications.  Let $\xi = c_1(L;\widetilde{\mathfrak M})$.  Then there is 
an isomorphism of
additive groups
$$
H^\dot_\cD (\mathbb P(E),\Z (p))_{\widetilde{\mathfrak M}}\cong \bigoplus_
{0\le j < r}
H^\dot_\cD (X,\Z (p))_{\mathfrak M}\cdot \xi^j,
$$
where $r$ is the rank of $E$; in other words, there is an injective map
$$
H^\dot_\cD (X,\Z (p))_{\mathfrak M}\to H^\dot_\cD (\mathbb P(E),\Z (p))_{
\widetilde{\mathfrak M}},
$$
and the right-hand side of the above is freely generated, as a module over 
the 
left-hand side, by $\{\xi^j: 0\le j < r\}$.
\end{proposition}
\begin{proof}
We want to reduce the assertion to more elementary cohomology for which one
already 
knows the splitting principle.  As in \ref{4.1.2}, there is a long exact 
sequence:
 \begin{multline}\label{4.1.2.A}
\dots\to H^{2p-1}(F^pA^\dot (\Xbar\log\,D))\to H^{2p-1}(X,\C /\Z (p))\cr
\to H^{2p}_\cD 
(X,\Z (p))_{\mathfrak M}\to H^{2p}(F^pA^\dot (\Xbar\log\,D))\to\dots
\end{multline}
Next, understanding $\xi$ to denote generically the Chern class of $L$ in any 
cohomology
group, 
multiply \ref{4.1.2.A} by $\xi^j$ ($0\le j < r$) and take the direct sum over
$j$ yields

\begin{equation}\label{4.1.3.A}
{\scriptsize
\setlength{\arraycolsep}{2pt}
\begin{array}{ccccc}
\oplus H^{2p-1-2j}(X,\C /\Z (p-j))\cdot\xi^j & \scriptstyle \to &  
\oplus H^{2p-2j}_\cD (X,\Z (p-j))_{\mathfrak M}
\cdot\xi^j & \scriptstyle \to & \oplus H^{2p-2j}(F^{p-j}A^\dot(\Xbar\log\,D))
\xi^j \cr
\downarrow && \downarrow && \downarrow\cr
H^{2p-1}(\mathbb P(E),\C /\Z (p)) & \scriptstyle \to & 
H^{2p}_\cD (\mathbb P(E),\Z (p))_{\widetilde{\mathfrak M}} & \scriptstyle\to &
H^{2p}(F^pA^\dot (\mathbb P(\Ebar)\log\,\psi^*D)) 
\end{array}
}
\end{equation}

Now, we know the splitting principle for $\Z$-coefficients, hence also for
$\Z (p)$- and
$\C$-coefficients, and then also for $(\C /\Z (p))$-coefficients.  In other 
words, 
the left vertical arrow in the display \ref{4.1.3.A} is an isomorphism. 
The easiest 
way to deal with the one on the right is to use Hodge theory: there is, in
general, a canonical surjection
\begin{equation}\label{hodge}
H^{2k}(F^{k}A^\dot (\Xbar\log\,D))\rightarrow F^{k}H^{2k}(X,\,\C ),
\end{equation}
and this 
is an isomorphism when $\Xbar$ is K\"ahler.  Thus, the right vertical arrow
is the
Hodge filtered version of the one for $\C$-coefficients, so it is an likewise
an isomorphism.
We conclude by the five-lemma that our assertion holds.
\end{proof}

The above determines \ref{CMBD}. We recall that $c_p(E;\mathfrak E)$ is, up 
to a sign, the 
coefficient in $\H^{2p}_\cD (X,\Z (p))_{\mathfrak M}$ of $\xi^{r-p}$ in the 
formula for 
$\xi^r$ in terms of the additive decomposition given in \ref{split}.  

Because of (\ref{6.1.4}), we obtain the generalization of Theorem \ref{thm2} 
to the 
K\"ahler case:
\begin{theorem}
Let $X$ be a quasi-c-K\"ahler manifold, $\mathfrak M$ an equivalence class of
K\"ahler 
compactifications of $X$, $E$ an $\mathfrak M$-extendable vector
bundle on $X$,
and $\nabla$ an $F^1$-connection on $E$.  Then $c_p^B (E;\mathfrak E)$ is the 
image of 
$\cs_p(E,\nabla)$ for all $p \geq 1$.  Here, $\mathfrak E$ is the unique 
equivalence 
class of $\mathfrak M$-extensions of $E$ implied by \ref{mero}. 
\end{theorem}

\begin{corollary}
The conclusion of Theorem 1 holds when $X$ is a quasi-c-K\"ahler manifolds
for the  $\mathfrak M$-DB  Chern classes of flat vector bundles that are
regular with respect to $\mathfrak E$.
\end{corollary}

\section{Proof of Theorem~\ref{thm4}}
\label{thm_3}
Throughout this section, the complexification $\g\otimes \C$ of a
real Lie algebra $\g$ will be denoted $\g_\C$. In addition, we make 
use of simplicial methods; all relevant definitions can be found in
Appendix~\ref{simp}.

\subsection{Continuous cohomology}\label{7.1} The continuous cohomology
$\Hcts\supdot(G,V)$ of a topological group $G$ with coefficients
in a real topological vector space $V$, on which $G$ acts, is
defined to be the cohomology of the complex $\Ccts\supdot(G,V)$,
whose elements in degree $m$ are $G$-equivariant continuous maps
$f: G^{m+1} \to V$ (with the usual coboundary map).

Now let $G$ be the real points of a reductive algebraic group defined
over $\R$, and $K$ a maximal compact subgroup.  We recall that the
space $G/K$ is contractible, and thus $G$ and $K$ are of the same
homotopy type. Let $\g = \k \oplus \p$ be the corresponding Cartan
decomposition of the Lie algebra of $G$. By van Est's Theorem~\cite{vE},
whenever $V$ is finite-dimensional with trivial $G$-action, there is
a canonical isomorphism
\begin{equation}\label{7.1.1}
\Hcts\supdot(G,V) \cong H\supdot(\g,K) \otimes V.
\end{equation}
Here, the right-hand side is the relative Lie algebra cohomology,
and is canonically isomorphic to
$\left(\wedge \supdot \p^\ast\right)^K \otimes V$, the complex
of $G$-invariant $V$-valued differentials forms on $G/K$.
The isomorphism can be realized as follows. Fix a point $e$ of $G/K$,
and endow the latter with a $G$-invariant Riemannian metric. For
$(g_0,\ldots,g_m)\in G^{m+1}$, define $\Delta_e(g_0,\ldots,g_m)$ to
be the ``geodesic simplex'' in $G/K$  with vertices $g_0e,\ldots,g_me$;
it is constructed inductively as the cone swept out by all geodesics
from $g_0 e$ to $\Delta_e(g_1,\ldots,g_m)$. (Note that, in general,
the order of the vertices makes a difference, since the curvature
may not be constant.) Define a homomorphism
$$
\wedge\supdot \p^\ast \to \Ccts\supdot(G,\R)
$$
by taking an $m$-form $\omega$ to the function
\begin{equation}\label{7.1.2}
(g_0,\ldots,g_m) \mapsto
\int_{\Delta_e(g_0,\ldots,g_m)} \omega.
\end{equation}
This induces \ref{7.1.1}.

We now recall the standard trick for computing the continuous
cohomology of a reductive group. The compact real form of
$\g\otimes \C$ is the Lie algebra
$\u = \k \, \oplus \, i\p \subseteq \g_\C$. Let $U$ be the
corresponding subgroup of
$G(\C)$. Since U is compact, there are canonical isomorphisms
$$
H\supdot(U/K,V) \cong H\supdot(\u,K;V) \cong
\left(\wedge\supdot(i\p)^\ast\right)^K \otimes V.
$$
Composing these isomorphisms with the isomorphism
$$
\wedge \supdot(i\p)^\ast \to \wedge \supdot \p^\ast
$$
induced by multiplication by $i$, we obtain an isomorphism
\begin{equation}\label{7.1.3}
H\supdot(U/K,\C) \mathop {\to}\limits ^{\sim}
\Hcts\supdot(G,\C),
\end{equation}
which carries $H^m(U/K,\R(p))$ onto the subspace $i^m\Hcts^m(G,\R(p))$
of $\Hcts^m(G,\C)$.

\subsection{The Borel regulator elements}\label{7.2} We now take $G$
to be $\GLn$, viewed as a real group, and we take $K$ to be $U(n)$. In
this case, one identifies the complexification of $G$ as
$\GLn \times \GLn$; the natural map $G \to G_\C$ is the homomorphism
$$
\iota : \GLn \to \GLn \times \GLn ,
$$
with $\iota (g) = (g,\gconj)$. The corresponding map on Lie algebras
\begin{equation}\label{7.2.1}
\gln \to \gln \oplus \gln
\end{equation}
takes $X$ to $(X,\Xconj)$; indeed, it is not hard to see that the
identification
\begin{equation}\label{7.2.2}
\gln \otimes _\R \C \cong \gln \oplus \gln
\end{equation}
is given by: for $X, Y \in \gln$,
\begin{equation}\label{7.2.3}
X\otimes 1 + Y\otimes i \mapsto \left(X + i Y, \overline X + i
\overline Y\right).
\end{equation}
It follows that the Lie algebra $\u_n$ of $U(n)$ appears in
\ref{7.2.1} as the ``anti-diagonal"
$$
\{ (X, -^tX)~|~ X \in \u _n \}.
$$

It is clear that the compact real form $U$ of $G_\C$ is
$U(n) \times U(n)$. The quotient space $U/K$ is therefore
$\left(U(n) \times U(n)\right)/ U(n)$, which we identify with
$U(n)$ via the map induced by inclusion as the first factor of the
product.  Since $\GLn$ and $U(n)$ are homotopy-equivalent, we obtain
from \ref{7.1.3} an isomorphism
$$
H\supdot(\GLn,\C) \cong \Hcts\supdot(\GLn,\C)
$$
which takes $H^m(\GLn,\R(p))$ onto $i^m\Hcts^m(\GLn,\R(p))$. This
yields a canonical isomorphism
\begin{equation}\label{7.2.4}
H^{2p-1}(\GLn,\R(p)) \mathop {\longrightarrow}\limits^{\sim}
\Hcts^{2p-1}(\GLn,\C/\R(p)).
\end{equation}

The real cohomology of $\GLn$ is an exterior algebra on generators
$x_1, \ldots, x_n$, where $x_p$ is in degree $2p -1$.  We can choose
$x_p$ canonically  by insisting that it be the unique (necessarily
primitive) class which transgresses, via
$$
d_{2p} : H^0(B\GLn ,H^{2p-1}(\GLn)) \to H^{2p}(B\GLn ),
$$
to the universal Chern class $c_p \in H^{2p}(B\GLn)$ in the Leray
spectral sequence associated to the universal $\GLn$-bundle; it
lies in
$$
W_{0}H^{2p-1}(\GLn,\Z(p)),
$$
as $\Z(p)$ is of weight $-2p$. Define
$b_p \in H^{2p-1}(B\GLn ^\delta ,\C/\R(p))$ to be the image of $x_p$
under the composite of \ref{7.2.4} and the natural map
$$
\Hcts\supdot(\GLn,\C/\R(p)) \to \Hcts\supdot(\GLn ^\delta ,\C/\R(p))
\cong H\supdot(B\GLn ^\delta ,\C/\R(p)).
$$
This is called the $p$-th {\it Borel regulator element}.\footnote{
Some prefer to define the Borel regulator element to be the class in
$H^{2p-1}(\GLn^\delta,\C/\R(p))$
that corresponds to the element $y_p$ of
$$
W_0 H^{2p-1}(\GLn,\Q(p)) \cong W_0 H^{2p-1}(GL_p(\C),\Q(p))
\cong  H^{2p-1}(\C^p - \{0\},\Q(p)) \cong \Q
$$
that takes the value 1 on the generator of $\pi_{2p-1}(GL_n(\C))$.
Since $x_p = \pm (p-1)! y_p$ \cite[(24.5.2)]{Sch}, it follows that $y_p$
corresponds, up to a sign, to twice the degree $p$ part of the Chern
character
$$
ch_p : K_{2p-1}(F) \to \Hdel^1(\Spec F, \R(p)).
$$
for $F$ a number field or $\R$ or $\C$.}

\subsection{The Weil algebra}\label{7.3} Let $G$ be a real reductive
group and $K$ any compact subgroup. Denote the corresponding Lie
algebras by $\g$ and $\k$. Consider the {\it (real) Weil algebra}
\cite{C1}
\begin{equation}\label{7.3.1}
W(\g) = \wedge \supdot \g^\ast \otimes S\supdot(\g^\ast),
\end{equation}
a d.g.a. that is a $\g$-module via the coadjoint action (here $S$ denotes
the symmetric algebra over $\R$).  Denote the set of invariant polynomials
$S\supdot(\g^\ast)^G$, which is a subalgebra of $W(\g)$, by $I\supdot(G)$.

Let $\omega \in A^1(P,\g)$ be a connection on a principal bundle
$P \to M$ over a (simplicial) manifold, which will be viewed as a map
$\g^\ast \to A^1(P)$. Then $\omega$ extends uniquely to a d.g.a.\
homomorphism
\begin{equation}\label{7.3.2}
k_P(\omega) : W(\g) \to A\supdot(P),
\end{equation}
such that for $\Phi \in S^p(\g^\ast)$,
$$
k_P(\omega)(\Phi) = \Phi(\Theta^p),
$$
where $\Theta$ denotes the curvature of the connection.

The {\it relative Weil algebra} $W(\g,K)$ \cite{C2} is the subspace of
``$K$-basic" elements
\begin{equation}\label{7.3.3}
\left[\wedge\supdot\left(\g/\k\right)^\ast \otimes
S\supdot(\g^\ast)\right]^K
\end{equation}
of $W(\g)$. One recovers \ref{7.3.1} when $K$ is trivial; $W(\g,K)$ is
a contravariant functor of pairs $(\g ,K)$, and it contains $I\supdot(G)$
for all $K$.  A connection $\omega \in  A^1(P,\g)$ on a principal bundle
$P \to M$ induces a d.g.a.\ homomorphism ({\it Chern-Weil homomorphism})
\begin{equation}\label{7.3.4}
k_P(\omega) : W(\g,K) \to A\supdot(P/K)
\end{equation}
by restriction of \ref{7.3.2}. The natural map $W(\g,K) \to I\supdot(K)$
(which factors through $I\supdot(G)$) is a quasi-isomorphism \cite{C2},
hence
$$
H\supdot (W(\g , K)) \cong H\supdot (K,\R).
$$
A consequence of this is that if $\Phi \in I\supdot(G)$ has the
property that its image in $I\supdot(K)$ is zero, then there exists
$T \in W(\g,K)$ such that $dT = \Phi$. Note that when $K$ is
{\it maximal} compact, $K \hookrightarrow G$ is a homotopy equivalence,
so $P/K \to P/G = M$ is a homotopy equivalence, and \ref{7.3.4} induces
a map
$$
H\supdot(G,\R) \to H\supdot(M,\R).
$$

For a real vector space $V$, denote the Weil algebra $W(\g)\otimes V$
with coefficients in $V$ by $W_V(\g)$, and the corresponding Weil algebra
in the relative case by $W_V(\g,K)$. In particular, we have the Weil
algebras $W_\C(\g)$ (quasi-isomorphic to $A\supdot (EG) \cong \C$) and
$W_\C(\g,K)$, as well as $W_{\R(p)}(\g,K)$.   When $\g$ and $\k$ are
{\it complex} Lie algebras, one can also define the complex Weil algebra
${\mathfrak W}(\g,K)$ by doing the linear algebra in \ref{7.3.3}
over $\C$ instead of $\R$.  We observe that when $\g$ is a real Lie
algebra, $W_\C(\g,K) \cong {\mathfrak W}(\g_\C,K_\C)$ as $\C$-algebras.

\subsection{The class $\tau_p$ and Cheeger-Simons classes}\label{7.4}
In this section, $\g$ is  $\gln$ viewed as a real Lie algebra, and
$\k = \u_n$. Let
$$
C_p(X) = (-1)^p ~\tr (\wedge^p X),
$$
the invariant polynomial (of degree p) on $\gln$ which determines the
$p$th Chern class. We can express it as
$$
C_p(X) = P_p(X) + Q_p(X),
$$
where
$$
P_p(X) = \frac{1}{2}\left[C_p(X) + (-1)^p C_p(\Xconj)\right]\text{ and }
Q_p(X) = \frac{1}{2}\left[C_p(X) - (-1)^p C_p(\Xconj)\right]
$$
are elements of $S^p(\g^\ast)$.  Note that $P_p$ is $\R (p)$-valued,
and $Q_p$ is $\R (p-1)$-valued. If $X\in \u_n$, then $\Xconj = -^t X$,
so $Q_p(X) = 0$.  It follows from the discussion in the previous section
that there exists $T_p \in W_{\R(p)}(\gln,U(n))$   such that $dT_p = Q_p$.

The connection on the universal flat bundle $P \to B\GLn ^\delta$ induces
a Chern-Weil homomorphism
$$
k : W_\C(\gln,U(n)) \to A\supdot\left(P/U(n)\right).
$$
Denote $k(T_p) \in A^{2p-1}(P/U(n))$ by  $\tau_p$; the curvature
$\Theta$ of $P$ is zero.  We have
$$
d\tau_p = Q_p(\Theta) = 0,
$$
and thus $\tau_p$ defines a class in
$$
H^{2p-1}(P/U(n), \C) \cong H^{2p-1}(\GLn^\delta,\C).
$$
The first task is to show:

\begin{proposition}\label{7.4.1}
The class of $\tau_p$ in $H^{2p-1}(\GLn^\delta,\C/\R(p))$ is the
negative of the universal real Cheeger-Simons class $\cs_p$.
\end{proposition}

\begin{proof}
Let $E\to B$ be the model of the universal $\GLn$-bundle constructed
from $P \to B\GLn^\delta$ (over the standard model of $B\GLn^\delta$),
as in \ref{6.1.5}. Recall that $E$ has a canonical connection, induced
by that of $P$, and there is a tautological map $B\GLn^\delta \to B$
which  classifies both the universal flat bundle and its connection.
Consider the Chern-Weil homomorphism
$$
k_E : W_\C(\gln,U(n)) \to A\supdot\left(E/U(n)\right)
$$
associated to $E$ and its connection, writing $\tau_p^E$ for $k_E(T_p)$.
Since the polynomial $P_p$ is $\R(p)$-valued,
$$
C_p(\Theta _E) \equiv Q_p(\Theta _E) \mod \R(p));
$$
and since $dT_p = Q_p$, it follows that $d\tau_p^E = Q_p(\Theta _E)$,
so $(Q_p(\Theta _E), -\tau_p^E)$ represents the ``universal
Cheeger-Simons class"
$$
\cs_p \in \widehat{H}^{2p}(E/U(n),\R(p)),
$$
i.e., that of the pullback of $E \to B$ to $E/U(n) \simeq B$.
The naturality of Chern-Weil construction implies that the diagram
$$
\begin{matrix}
W_\C(\gln,U(n)) & \mathop {\to}\limits^{k_E} &
A\supdot\left(E/U(n)\right)\cr
& \mathop {\searrow}\limits_k & \downarrow \cr
& & A\supdot\left(P/U(n)\right) \cr
\end{matrix}
$$
commutes. It follows that the Cheeger-Simons class of the pullback
of $P \to B\GLn^\delta$ to $P/U(n)$ is represented by
$(Q_p(\Theta ),-\tau_p) = (0,-\tau_p)$, and this gives the desired
assertion.
\end{proof}

\subsection{$\tau_p$ and $T_p$}\label{7.5}
The next step is to obtain a formula for the cocycle given by
$\tau_p$ in terms of $T_p$. For this we will need an explicit section
of the bundle $P/U(n) \to B\GLn^\delta$. It is more convenient to write
down a section of $E\GLn/U(n) \to B\GLn$ and pull it back to the
universal flat bundle.  Actually, we will write down a section of the
bundle $\overline \pi : EG/K \to BG$, where $K$ is a maximal compact
subgroup of an arbitrary real reductive group $G$, following \cite{D1}.

The $m$-simplex $G^{m+1}$ of the usual model of $EG$ (see \ref{6.1.5},
and take $X\subdot$ to be a point) is identified with $G^{m+1}$ in
the {\it inhomogeneous}, or {\it reduced}, model via the map
$$
(g_0,g_1,\ldots,g_m) \mapsto (g_0,g_1g_0^{-1},\ldots,g_mg_{m-1}^{-1}).
$$
This map is $G$-equivariant with respect to the diagonal right
$G$-action on the left-hand side, and the right action of $G$ on the
first factor of the right-hand side. It follows that in the reduced
model, the space of $m$-simplices of $EG/K$ is $G/K \times G^m$. In
particular, when $G=K$, we have that the space of $m$-simplices of
the reduced model of $BG$ is $G^m$.

Suppose now that $G$ is reductive and that $K$ is a maximal compact
subgroup. For each $m$, define a map (which is a homotopy equivalence)
\begin{equation}\label{7.5.1}
s_m : G^m \times \Delta^m \to G/K \times G^m \times \Delta^m
\end{equation}
by
$$
\scriptstyle
\left((g_1,\ldots, g_m),(t_0,\ldots,t_m)\right) \mapsto
\left((\widetilde \Delta_e ((g_1,g_2g_1,\ldots,g_mg_{m-1}\cdots
g_1),(t_0,\ldots,t_m)),
(g_1,\ldots, g_m),(t_0,\ldots,t_m)\right),
$$
where $\widetilde \Delta_e : G^m \times \Delta^m \to G/K$ parametrizes
the geodesic simplices defined in \ref{7.1.2}.

\begin{proposition}{\cite[p.~241]{D1}}\label{7.5.2}
The maps \ref{7.5.1} are compatible with the face maps and
therefore induce a section $s: BG \to (EG)/K$ of $\overline \pi$.
\endproof
\end{proposition}

We now show that the cohomology class on $\GLn^\delta$ defined by
the continuous cohomology class corresponding to
$T_p \in \wedge^{2p-1}\p^\ast$ equals the cohomology class defined by
$s^\ast\tau_p$. With respect to the inhomogeneous model
(see \ref{7.5.1}), the connection form of the universal flat bundle
is given locally by the Maurer-Cartan form in the vertical direction,
and is zero in the horizontal directions. Consequently, with respect
to the inhomogeneous model, the restiction of
$\tau_p$ to $\GLn/U(n)\times G^m \times \Delta^m$ satisfies
$$
\tau_p \in A^{2p-1}(\GLn/U(n)) \subseteq
A^{2p-1}(\GLn/U(n)\times G^m \times \Delta^m),
$$
and that it is the image of $T_p$ under the natural inclusion
$$
\wedge^{2p-1} \p^\ast \hookrightarrow A^{2p-1}(\GLn/U(n)).
$$
It follows immediately that the value of the cocycle
$$
s^\ast\tau_p\in A^{2p-1}(B\GLn^\delta)
$$
on the simplex of $B\GLn$ with vertices $(g_0,\ldots,g_m)$ is given
by the integral
$$
\int_{\Delta_e(g_0,\ldots,g_m)} T_p.
$$

\subsection{$T_p$ and the Borel regulator element}\label{7.6}
In the final step, the idea is to use the complexification to show
that $T_p$ can be interpreted as a class that transgresses to the
$p$th Chern class of $E \to B$.

An invariant polynomial on $\gln$ defines, in an obvious way, an
invariant function on $\g_\C$.  Using \ref{7.2.2} to write the
elements of $\g_\C$ as pairs $(X,Y)$, with $X,Y \in \gln$, we have
for $Q_p$ of (7.4):
\begin{equation}\label{7.6.1}
Q_p(X,Y) = \frac{1}{2}\left[C_p(X) - (-1)^p C_p(Y)\right].
\end{equation}
{}From this, we see that for all $X \in \gln$,
\begin{equation}\label{7.6.2}
C_p(X) = 2\,Q_p(X,0).
\end{equation}
Likewise, we can view $T_p$ (also from \ref{7.4}) as an element
of the complex Chern-Weil algebra ${\mathfrak W}(\g_\C , U(n)_\C )$,
which then determines elements of $W_\C(\u_n\oplus\u_n,U(n))$ and
$$
W_\C(\u_n \oplus \{ 0\} ) \cong W_\C(\u_n),
$$
which we continue to denote $T_p$.

{}From \ref{7.6.2}, it follows that, in $W_\C(\u_n)$, the relation
$2\,dT_p = C_p$ holds. In other words, the invariant form on $U(n)$
defined by $2\,T_p$ represents a class which transgresses to the
$p$th Chern class. The corresponding class in
$$
\wedge \supdot \left(\u_n\right)^\ast \mathop {\to} \limits^\sim
_{{pr_1}^\ast} \wedge \supdot \left((\u_n \oplus \u_n)/\u_n\right)^\ast
\cong \wedge \supdot (i\p^\ast)
$$
is the image of $2\,T_p$ in $\wedge^{2p-1}(i\p^\ast)$.
By the definition given in \ref{7.2}, the Borel element is
represented by the differential form
$$
i^{2p-1}\left(2\,T_p|_{\wedge^{2p-1}i\p}\right)
= (-1)^{2p-1}\left(2\,T_p|_{\wedge^{2p-1}\p}\right) = -2\,T_p \in
\wedge^{2p-1} \p^\ast.
$$

Combining the above with \ref{7.4} and \ref{7.5}, we have shown
that the Borel element is twice the Cheeger-Simons class, and
Theorem~3 is proved.

\section{Towards a Proof of Conjecture~4}
\label{discussion}

In this section we make use of simplicial methods.  Again, see 
Appendix~\ref{simp} for definitions and constructions.

\subsection{The goal}\label{8.1}
We would have preferred to have 
Theorem~\ref{thm2} and Conjecture~\ref{conj}
as special cases of the simplicial
analogue of Theorem~\ref{thm2}, which we have not been able to prove:

\begin{statement}\label{6.1.1}
Suppose that $ \pi : \Pbar_\dot \rightarrow \Xbar_\dot$ is a morphism
in the category of smooth simplicial varieties, with $\Xbar_\dot$
complete, and that  $D\subdot$ and $Q\subdot = \pi^{-1}(D\subdot)$
and are divisors with normal crossings in
$\Xbar_\dot$ and $\Pbar_\dot$, respectively. Let
$P\subdot = \Pbar_\dot- Q\subdot$ and
$X\subdot = \Xbar_\dot - D\subdot$. If the restriction
$P\subdot \rightarrow X\subdot$ of $\pi$ is a
principal $GL_n(\C)$-bundle with $F^1$-connection ---
i.e., the connection form $\omega$ satisfies
$$
\omega \in F^1 E^1(|\Pbar_\dot|\log Q_\dot,\End(\Ebar_\dot)),
$$
then the image of its Cheeger-Simons class
$\cs_p$ in $\Hdel^{2p}(X\subdot ,\Z(p))$ under the natural homomorphism
$$
F^p \widehat H^{2p} (\Xbar_\dot\log D,\Z (p))\longrightarrow
\Hdel^{2p}(X\subdot, \Z(p))
$$
is the Beilinson Chern class $c^B_p$.
\end{statement}

\begin{remark}\label{6.1.2}
Note that we are {\it not} assuming in the above that $\Pbar_\dot$ is
proper over $\Xbar_\dot$, but only that it is so ``in the horizontal
direction".  In practice, $P\subdot$ would be the frame
bundle of a vector bundle $E\subdot$ on $X\subdot$, and $\Pbar_\dot$
the frame bundle of an extension $\Ebar\subdot$ of $E\subdot$ to
$\Xbar_\dot$.
\end{remark}

The following lemma, whose proof is quite direct, provides good evidence for
the conjecture.
(We retain the previous notation.)

\begin{lemma}\label{6.1.3}
Suppose that one has two connections $\nabla_0$ and $\nabla_1$
on $E\subdot$, and denote by $\Theta_0$ and $\Theta_1$ the
respective curvature forms.  Write $\omega = \nabla_1 - \nabla_0$
and $\eta_p$ for the solution of
$$
d\eta_p = c_p(E\subdot, \nabla_1) - c_p(E\subdot, \nabla_0)
$$
given by \ref{3.5.5}. If
$$
\Theta_0, \Theta_1 \in F^1A^2(|\Xbar_\dot|\log D\subdot,
\End(\Ebar\subdot))\,\, {\it and }\,\,
\omega \in F^1(A^1(|\Xbar_\dot|\log D\subdot, \End(\Ebar\subdot))),
$$
then
\begin{enumerate}
\item[(i)]$\eta_p \in F^p A^{2p-1}(|\Xbar_\dot|\log D\subdot))$;
\item[(ii)]$c_p(E\subdot, \nabla_0)$ and $c_p(E\subdot,\nabla_1)$
represent the  same cohomology class in
$$
H^{2p}(F^pA^{\dot}(|\Xbar_\dot|\log D\subdot))=F^pH^{2p}(|X\subdot|,\C).
$$
\end{enumerate}
\end{lemma}

\begin{proof} This is a direct consequence of \ref{3.5.4} and
\ref{3.5.5}, as the assumptions imply that $\Theta_t$ is in $F^1$.
\end{proof}

The following ``refinement" of \ref{6.1.3}, which we will soon
prove, is a necessary condition for \ref{6.1.1}:

\begin{proposition}\label{6.1.4} Under the conditions of
\ref{6.1.3}, these two Cheeger-Simons classes have the same
image in $\Hdel^{2p}(X\subdot,\Z(p))$.
\end{proposition}

The basic theme in any proof of \ref{6.1.1} is to reduce to the
universal case (where the assertion is a tautology; see
\ref{3.2.4}, and \ref{6.1.10} below), invoking functoriality.
We wanted to make use of the following device:

\begin{proposition}\label{6.1.5}
If $X\subdot$ is a simplicial manifold and
$P\subdot \rightarrow X\subdot$ is a principal
$GL_n(\C)$-bundle with connection given by
$$
\omega \in A^1(|P\subdot|, \Ad(\gln)),
$$
then there exists a smooth bisimplicial variety $B\subdotdot$ which
has the homotopy type of $BGL_n(\C)$ and a $GL_n(\C)$-principal bundle
$U\subdotdot \rightarrow B\subdotdot$ with connection
$$
\omega_U \in A^1(|U\subdotdot|,\Ad(\gln))
$$
and a morphism of $GL_n(\C)$-bundles
$$
\begin{matrix}
P\subdot & \buildrel {G} \over \longrightarrow & U\subdotdot \cr
               \downarrow &   &   \downarrow \cr
X\subdot  & \buildrel {g} \over \longrightarrow &  B\subdotdot
\end{matrix}
$$
such that  $G^{\ast}\omega_U = \omega$. Moreover, if $X_\dot$ is a smooth
simplicial variety such that
$$
\omega \in F^1A^1(|P\subdot|,\Ad(\gln)),
$$
then we can choose $\omega_U$ to satisfy
$$
\omega_U \in A^{1,0}(|U\subdotdot|,\Ad(\gln)).
$$
\end{proposition}

\begin{proof} Consider the bisimplicial variety $U\subdotdot$
whose simplicial  variety of $m$ simplices is the simplicial variety
$(P\subdot)^{m+1}$ with face maps
\begin{align*}
d_i : (P\subdot)^{m+1} & \rightarrow (P\subdot)^m
\qquad 0\le i \le m\cr
(u_0, \ldots, u_m) & \mapsto (u_0, \ldots, \hat{u}_i,
\ldots, u_m).
\end{align*}
The geometric realization of this simplicial variety is contractible
as it is the coskeleton of the trivial covering (cf.\ \cite{AM} or
\cite[p.~107]{Se}).
\footnote{This may be proved directly as follows. First show that the
geometric realization of the simplicial space is simply connected. This is
not difficult, as the fundamental group depends only on the 2-skeleton.
One then shows that this space has trivial integral homology, and is
therefore contractible, by standard arguments.}  The free $GL_n(\C)$
action on $P\subdot$ induces a free $GL_n(\C)$ action on $U\subdotdot$.
Let $B\subdotdot$ be the quotient of $U\subdotdot$. The resulting
principal bundle $U\subdotdot \rightarrow B\subdotdot$ is a model of
the universal $GL_n(\C)$-bundle.

Define a connection $\omega_U \in A^1(|U\subdotdot|,\gln)$ on the
bundle $|U\subdotdot| \rightarrow |B\subdotdot|$ by the compatible
family of 1-forms
$$
\sum_{j=0}^m t_j\pi^{\ast}_j \omega \in A^1(|P\subdot|^{m+1} \times
\Delta^m),
$$
where $(t_0, \ldots, t_m)$ are the barycentric coordinates of
$\Delta^m$, and $\pi_j : |P\subdot|^{m+1} \rightarrow |P\subdot|$
denote the canonical projections, $0 \leq j \leq m$.

The canonical isomorphism $P\subdot \rightarrow U_{\dot 0}$
induces a $GL_n(\C)$-equivariant map
$G: P\subdot \rightarrow U\subdotdot$, and therefore a map
$g: X\subdot \rightarrow B\subdotdot$ such that the
pair of maps $(g,\, G)$ classifies the bundle and the connection.
The remaining assertion is easy to verify.
\end{proof}

\begin{remark}\label{6.1.7}
In other words, given a principal $GL_n(\C)$-bundle with connection
over a complete smooth simplicial variety , we can build a model of the
universal $GL_n(\C)$-bundle, and put a connection on this bundle, such that
a {\it fixed} mapping of our variety into this model of $BGL_n(\C)$
simultaneously classifies the bundle and its connection.  Unfortunately the 
above does not yield a proof of \ref{6.1.1}, for the 
Maurer-Cartan 
forms of $GL_n(\C)$, and hence the connection constructed in \ref{6.1.5}, are 
not   
logarithmic at infinity for $n>1$. This was our fundamental obstacle.
\end{remark}

However, \ref{6.1.5} does yield at once the following useful observation:

\begin{proposition}\label{6.1.6}
Suppose that $\nabla_1$ is a second connection on $E\subdot$ and that
$\cs_p (E\subdot,\nabla_0)$ is represented by
$(c_p (E\subdot, \nabla_0), -y)$, as in \ref{3.5.17}.  Let $\eta_p$
be as in \ref{3.5.5}. Then
$\cs_p (E\subdot, \nabla_1)$ is represented by
$$
(c_p (E\subdot, \nabla_1), -(y + \iota_{\Z(p)}(\eta_p ))) =
(c_p (E\subdot, \nabla_0), -y) + (d\eta_p, -\iota_{\Z(p)}(\eta_p ))).
$$
\end{proposition}

\begin{proof} By the functoriality of $c_p$, $\eta_p$ and
$\cs_p$, it suffices to check this ``universally", on the model of
$BGL_n(\C)$ from \ref{6.1.5}.  But there, the Cheeger-Simons class
is completely determined by the Chern form (see \ref{3.4}), so it
comes down to the calculation:
$$
\delta (y + \iota_{\Z(p)}(\eta_p)) = \iota_{\Z(p)}(c_p
(E\subdot , \nabla_0)) +
\iota_{\Z(p)}(d \eta_p) = \iota_{\Z(p)}(c_p (E\subdot ,
\nabla_1)).
$$
\end{proof}

Using \ref{6.1.6}, we can now give:

\begin{proof}[Proof of \ref{6.1.4}]  In the DB-complex of $X\subdot$, we
have that:
$$
(c_p (E, \nabla_0), -y)
 - (c_p (E, \nabla_1), -(y + \iota_{\Z(p)}(\eta_p ))) =
(-d\eta_p , \iota_{\Z(p)}(\eta_p )) = D(\eta_p , 0),
$$
so \ref{6.1.4} follows.
\end{proof}

The simplicial version of \ref{4.2.1} is
\begin{equation}\label{6.1.8}
F^p\widehat H^k(|\Xbar_\dot|\log D_{\dot}, \Lambda )
= H^{k-1}(\cone\{ F^pA^{\geq k}(|\Xbar_\dot|\log \Dbar_\dot) \rightarrow
S^{\dot}(|X\subdot |, \C / \Lambda )\}).
\end{equation}
Taking $k = 2p$, and $\Lambda = \Z(p)$, and recalling the definitions in
Section~\ref{db_coho}, we obtain the following diagram, in which the rows
are exact:
\begin{equation}\label{6.1.9}
{\setlength{\arraycolsep}{4pt}
\begin{matrix}
\scriptstyle 0 &  \scriptstyle \to &
\scriptstyle H^{2p-1}(|X\subdot|,\C/\Z(p)) &
\scriptstyle \to &
\scriptstyle F^p \widehat H^{2p}(\Xbar_\dot\log D_{\dot}, \Z(p) ) &
\scriptstyle \to & \scriptstyle
F^pA^{2p}(|\Xbar_\dot|\log D\subdot,\Z(p))_{cl}
& \scriptstyle \to & \scriptstyle 0 \cr
& & || & & \scriptstyle \downarrow  & & \scriptstyle \downarrow &  & \cr
& & \scriptstyle H^{2p-1}(|X \subdot|,\C /\Z(p) ) & \scriptstyle \to
&\scriptstyle \Hdel^{2p}(X \subdot , \Z(p) ) & \scriptstyle \to &
\scriptstyle F^pH^{2p}(|X\subdot|,\C) &  &
\end{matrix}
}
\end{equation}
We note that the following is easily read from \ref{6.1.9}:

\begin{proposition}\label{6.1.10} If $H^{2p-1}(|X \subdot |,\C /\Z(p))=0$,
then the assertion in \ref{6.1.1} holds for every holomorphic vector
bundle on $|X \subdot |$.
\end{proposition}

Next, consider the diagram
$$
\begin{matrix}
\Hdel^{2p}(\Xbar_\dot, \Z(p)) & \longrightarrow & \Hdel^{2p}(X\subdot,
\Z(p)) \cr
\uparrow & & \uparrow \cr
 F^p \widehat H^{2p} (\Xbar_\dot, \Z (p)) & \longrightarrow &
 F^p \widehat H^{2p} (\Xbar_\dot \log D_{\dot}, \Z (p))
\end{matrix}
$$
Here, we would like to use \ref{6.1.4} to replace the given connection by
one that extends without singularity to the compactification, without
changing the image in Deligne-Beilinson cohomology; the desired
conclusion then follows by functoriality.  This can always be done
when $X\subdot$ is just a single algebraic manifold $X = \Xbar - D$.
Thus, the proof of Theorem~\ref{thm2} (the non-simplicial version of
\ref{6.1.1}) reduced to the compact case (see \ref{6.2}).

\begin{remark}\label{6.1.11}
It should be apparent that the preceding development goes through
verbatim to the setting of simplicial complex manifolds $X \subdot$
that are of the form $\Xbar_\dot - D\subdot$, with $\Xbar_\dot$
compact K\"{a}hler and $D\subdot$ a divisor with normal crossings in
$\Xbar_\dot$.  In fact, even the condition that $\Xbar_\dot$ be K\"{a}hler
can be dispensed with; one simply makes the distinction between the two
sides (isomorphic when $\Xbar_\dot$ is K\"{a}hler) of
$$
F^pH^{\dot} (X\subdot, \C ) \longleftarrow H^{\dot}
(F^pA^\dot (|\Xbar_\dot |\log D\subdot ))
$$
and
$$
H^{\dot} (X\subdot, \C )/F^pH^{\dot} (X\subdot, \C )
\hookrightarrow
H^{\dot} (A^\dot (|\Xbar_\dot |\log D\subdot )/F^p),
$$
and takes the right-hand member wherever we have written the left.
\end{remark}

We present in Appendix D a couple of plausible techniques we tried in our
unsuccessful attempt to prove \ref{6.1.1}.

} 

\appendix
\section{Fiber integration and the homotopy formula}
\label{3A}

 In this appendix, we present some basic facts needed in 
Section~\ref{cs}.
Let $p : I \times M \to M$ denote the canonical projection of the
product of the unit interval with an $m$-manifold $M$ onto $M$. For
all $i$, there is a linear mapping
\begin{equation}\label{3.A.1}
B : A^{i + 1} (I \times M) \to A^i (M)
\end{equation}
characterized by the formula
\begin{equation}\label{3.A.2}
T (B \varphi) = ([I] \times T)
\end{equation}
for all compactly supported $i$-currents $T$ on $M$.  Whenever $T$ is
the current $T_\psi$ defined by
\begin{equation}\label{3.A.3}
T_\psi (\varphi) = \int_M \varphi \wedge \psi \qquad (\psi \in
A_c^{m-i}(M)),
\end{equation}
for which
\begin{equation}\label{3.A.4}
(i)\qquad \partial T_\psi = (-1)^{i+1} T_{d\psi}
\text{ and }\quad (ii)\qquad [I] \times T_\psi = T_{p^\ast \psi},
\end{equation}
\ref{3.A.2} becomes $\int_M B \varphi \wedge \psi = \int_{I \times M}
\varphi \wedge p^\ast \psi$.  In fact, in general, $B \varphi$ is just
the integral over $I$ of $(\partial_t \rfloor \varphi)$.  The mapping
$B$ is functorial, in the sense that if $g: N \to M$ is smooth, inducing
$\widetilde g : I \times N \to I \times M$, then
\begin{equation}\label{3.A.5}
B (\widetilde g^\ast \varphi) = g^\ast B\varphi.
\end{equation}

The identity
\begin{equation}\label{3.A.6}
\partial ([I] \times T) = \partial [I] \times T - [I] \times
\partial T\
\end{equation}
yields at once
\begin{equation}\label{3.A.7}
B (d \varphi) + d (B \varphi) = \varphi_1 - \varphi_0,
\end{equation}
where $\varphi_t = \varphi |_{\{ t \} \times M}$.
Applied to $\varphi = h^\ast \omega$, where $h : M \times I \to N$
is a smooth homotopy, and $\omega \in A^i (N)$, \ref{3.A.7} gives
the {\it homotopy formula}:
\begin{equation}\label{3.A.8}
(B \circ h^\ast) d\omega + d (B \circ h^\ast) \omega = h_1^\ast
\omega - h_0^\ast \omega.
\end{equation}
The dual version for currents is
\begin{equation}\label{3.A.9}
\partial (h_\ast ([I] \times T)) + h_\ast ([I] \times \partial T) =
(h)_\ast T- (h_0)_\ast T.
\end{equation}

There is a parallel treatment for singular cochains with any
coefficients, as follows.  A smooth singular simplex
$\sigma : \Delta^i \to M$ defines a current:
\begin{equation}\label{3.A.10}
[\sigma] = \sigma_\ast [\Delta^i],
\end{equation}
such that $[\partial \sigma] = \partial [\sigma]$.
Using the standard triangulation of $I \times \Delta^i$, a smooth
mapping
$$
\tau: I \times \Delta^i \to N
$$
defines a singular chain $[\tau]$ on $N$. Taking $\tau$ of the form
$$
h_\ast (1 \times \sigma): I \times \Delta^i \to I \times M \to N,
$$
we have
\begin{equation}\label{3.A.11}
\partial h_\ast (1 \times \sigma) + h_\ast (1 \times \partial
\sigma) = (h_1)_\ast \sigma - (h_0)_\ast \sigma,
\end{equation}
which under \ref{3.A.10} recovers \ref{3.A.9}.  Define, for a smooth
singular cochain $f$ on $I \times M$,
\begin{equation}\label{3.A.12}
(Bf) (\sigma) = f (1 \times \sigma)
\end{equation}
(cf.\  \ref{3.A.2}); putting $f = h^\ast \tau$, we get from
\ref{3.A.11}
\begin{equation}\label{3.A.13}
B (h^\ast \delta \tau) + \delta (Bh^\ast \tau) = h_1^\ast \tau -
h_0^\ast \tau.
\end{equation}
When $\tau$ is a differential form (see \ref{3.1}), one recovers
\ref{3.A.8}, i.e.:

\begin{proposition}\label{3.A.14} The homotopy formula for smooth
singular cochains, when restricted to differential forms, coincides
with the homotopy formula for differential forms. \endproof
\end{proposition}

\section{Universal weak $F^1$-connections}
\label{4A}
We show that if one ceases to impose conditions at infinity, it is possible
to construct
a universal ``$F^1$-connection.''

We start with a definition:
\begin{definition}\label{4.2.w} 
A {\it weak $F^1$-connection} on $E$ is a connection for 
which the local connection forms and therefore also the curvature lie in
$F^1 A^\dot (X,\End(E))$.
\end{definition}

\noindent Imposing the growth condition along $D$, as in \ref{4.2.1}, in 
\ref{4.2.w}
gives \ref{4.F}.

We introduce the notion of a {\it universal weak $F^1$-connection}.
The discussion proceeds as in that of \ref{2.2}.
In the setting of \ref{2.2.1}, suppose that $Y$ is a {\it complex}
manifold, and take $G = GL_n(\C)$. Let $P \to Y$ be a holomorphic
principal $GL_n(\C)$-bundle, and
$$
p: T^{1,0}P \to T^{1,0}Y
$$
be the induced map of holomorphic tangent bundles.  We put
\begin{equation}\label{5.6.1}
{\widetilde Y}^{(1)} = \{\widehat \alpha \in\Hom(T^{1,0}Y,T^{1,0}P/G):
\widehat  p\circ
\widehat \alpha = \id_{T^{1,0}Y}\},
\end{equation}
where ${\widehat p}: T^{1,0}P/G \to T^{1,0}Y$ is induced by $p$. Then
${\widetilde Y}^{(1)}$ can be identified with a submanifold of the
space ${\widetilde Y}$ in \ref{2.2.1} in the following way.  One makes
use of the canonical splitting $(TY)_{\C} = T^{1,0}Y \oplus T^{0,1}Y$.
The projection $q^{(1)}:{\widetilde Y}^{(1)}\to Y$ then becomes a
complex affine space sub-bundle of $q:\widetilde Y \to Y$ (via
$\alpha = \widehat \alpha \oplus 0$); a local holomorphic
trivialization of $P$ gives rise to the local description
(cf.~\ref{2.2.6}):
\begin{equation}\label{5.6.2}
{\widetilde Y}^{(1)} \cong \Hom(T^{1,0}Y,\g).
\end{equation}
There is a tautological {\it holomorphic} connection on ${q^{(1)*}}P$.

\begin{proposition}\label{5.6.3}
The weak $F^1$-connections on $P$ are in one-to-one
correspondence with the $C^{\infty}$ sections of $q^{(1)}$; the
holomorphic sections of $q^{(1)}$ correspond to the holomorphic
connections on $P$.
\end{proposition}

One must remember that the notion of weak $F^1$-connection is respected
under pullback only by {\it holomorphic} maps. For a holomorphic map
$f:X\to Y$, the pullback connection on $f^*P$ is induced by
$$
\Hom(T_{f(x)}^{1,0}Y,\g) \to \Hom(T_x^{1,0}X,\g).
$$
Taking $Y = G_{\C}(n)$, we deduce (cf.~\ref{2.2.9}, \ref{2.2.10}):

\begin{corollary}\label{5.6.4}
If a holomorphic vector bundle $E$ is classified by a
holomorphic immersion $f: X \to G_{\C}(n)$, then every weak $F^1$-connection
on $E$ is the pullback of the tautological weak $F^1$-connection on
${\widetilde G}_{\C}(n)^{(1)}$ by a $C^{\infty}$ lifting
${\widetilde f}: X \to {\widetilde G}_{\C}(n)^{(1)}$ of $f$.
\end{corollary}

\section{Simplicial de~Rham Theory}
\label{simp}

In this appendix we give a resume of simplicial manifolds, tailored
to our present needs.  The principal reference is \cite{D2}.

\subsection{Simplicial manifolds}\label{A.1}
Let $\simp$,be the category whose objects are the finite ordinals
$$
[n] = \{ 0 < 1 < \ldots < n \} \quad n \in \N
$$
and whose morphisms are order-preserving {\it monomorphisms}. Of
special importance are the face maps $d_j : [n-1] \to [n]$, the
unique morphism that omits the value $j$.

View $[n]$ as the set of the set of vertices of the standard
$n$-simplex $\Delta^n$.  Then, each morphism $f : [n] \to [m]$ of
$\simp$ induces a simplicial map
$$
| f | : \Delta^n \to \Delta^m.
$$

\begin{definition}\label{A.1.1} A {\it (strict) simplicial}
(resp.\ {\it cosimplicial}) object of a category $\cC$ is a
contravariant (resp.\  covariant) functor $\simp$  $\to \cC$.
Morphisms between strict (co)simplicial objects in $\cC$ are
natural transformations of functors.

\end{definition}
Denote the category of smooth manifolds, each of whose connected
components are second countable, and smooth maps by $\cM$.

\begin{definition}\label{A.1.2} A {\it simplicial manifold} $M_\dot$
is a strict simplicial object of $\cM$. The object of $\cM$
corresponding to the ordinal $[n]$ will be denoted by $M_n$. The smooth
map induced by the morphism $f : [n] \to [m]$ will be denoted by
$$
f_M : M_m \to M_n.
$$
\end{definition}

The {\it geometric realization} $| M_\dot |$ of the
simplicial manifold $M_\dot$ is defined to be the topological space
$$
\coprod_{n \geq 0} (M_n \times \Delta^n)/\sim
$$
endowed with the quotient topology, where $\sim$ is the equivalence
relation which identifies the points
$$
(x, | f | \xi) \in M_m \times \Delta^m \quad \text { and }
\quad (f_M (x), \xi) \in M_n \times \Delta^n
$$
whenever $f : [n] \to [m]$ is a morphism.

Observe that simplicial sets can be viewed as simplicial manifolds.
In particular, the standard simplicial model of the classifying space
of  $GL_n(\C)^\delta$, the general linear group {\it with the discrete
topology}, is a simplicial manifold.

There  is a fully faithful imbedding of $\cM$ into $\cM^\simp$, the
category of simplicial manifolds: The manifold $M$ is taken to the
simplicial manifold $M_\dot$ with
$$
M_n =
\begin{cases}
M & n = 0 \cr  \emptyset & n > 0.
\end{cases}
$$

\subsection{The de~Rham functor}\label{A.2} Suppose that $C^\dot$
is a strict cosimplicial cochain complex of $\C$-vector spaces.
Denote the cochain complex corresponding to $[n]$  by $C^\dot [n]$.

\begin{definition}\label{A.2.1}{\cite[(5.3)]{H}}
The {\it de~Rham complex} $DC^\dot$ of $C^\dot$ is the total
complex of the double complex $D^{\dot\dot}C^\dot$, where
$D^{s, t} C^\dot$ consists of those elements $(w_n)$ of
$$
\prod_{n \geq 0} A^s (\Delta^n) \otimes C^t [n]
$$
that satisfy the compatibility condition:
$$
(| f |^\ast \otimes \id) w_n = (id \otimes f_\ast) w_m \in
A^s(\Delta^m) \otimes C^t [n],
$$
whenever $f : [m] \to [n]$ is a morphism.  The differentials in the
double complex arise from exterior differentiation in the simplicial
direction $s$ and from the differential of $C^\dot$ in the $t$
direction.
\end{definition}

Now suppose that $M_\dot$ is a simplicial manifold.  Applying the
usual de~Rham complex functor to $M_\dot$ we obtain a cosimplicial
cochain complex $A^\dot (M_\dot)$ whose value on the ordinal
$[n]$ is $A^\dot (M_n)$.

\begin{definition}\label{A.2.2} The {\it de~Rham complex
$A^\dot| M_\dot |$ of the simplicial manifold} $M_\dot$
is the de~Rham complex $DA^\dot (M_\dot)$ of
$A^\dot (M_\dot)$.  It is naturally bigraded:
$$
A^k | M_\dot | =
\bigoplus_{s + t =k } A^{s, t} |  M_\dot|,
$$
where
$$
A^{s, t} | M_\dot | = D^s A^t (M_\dot).
$$
\end{definition}

If $M_\dot$ is a simplicial set, then $A^\dot|M_\dot|$
is just the Thom-Whitney de~Rham complex of $\C$-valued forms
on $M_\dot$ (cf.\ \cite{Sw}, \cite{W}).

\begin{proposition}\label{A.2.3}
The association of $A^\dot |M_\dot|$ to the simplicial manifold
$M_\dot$ defines a contravariant functor from the category of
simplicial manifolds into the category of commutative differential
graded algebras. \endproof
\end{proposition}

\subsection{The de~Rham Theorem}\label{A.3} In this section we sketch
the de~Rham theorem for simplicial manifolds.  It is first convenient
to introduce the functor $C^\dot$ from the category of (strict)
cosimplicial abelian groups into the category of cochain complexes:
The cosimplicial abelian group $A$ is taken to the cochain complex
$C^\dot A$, where
$$
C^n A = A [n]
$$
and the differential $\delta : C^n A \to C^{n+1} A$ is the
alternating sum
$$
\delta = \sum_{j=0}^{n+1} (-1)^j d^j
$$
of the coface maps $d^i : A [n] \to A [n+1]$.
When $A^\dot$ is a strict cosimplicial cochain complex, $C^\dot
A^\dot$ is a double complex.

The classical de~Rham theorem states that if M is a smooth
manifold, the integration map
$$
\int : A^\dot (M) \to S^\dot (M, \C)
$$
is a quasi-isomorphism which induces an algebra isomorphism
on cohomology.

Now, to each simplicial manifold $M_\dot$,
we can associate the strict cosimplicial cochain complex $S^\dot
(M_\dot)$ whose value on the ordinal $[n]$ is $S^\dot (M_n)$.
Applying the
cochain functor $C$ above, we obtain a double complex $C^\dot
S^\dot (M_\dot)$.  Denote the associated single complex (total
complex)
by $S^\dot  (M_\dot)$.
Note that if $M_\dot$ is a strict simplicial set, then $S^\dot
(M_\dot)$ denotes the usual complex of simplicial cochains
associated to
$M_\dot$.

For each simplicial manifold $M_\dot$, define a $\C$-linear
function
\begin{equation}\label{A.3.2}
A^{s, t} | M_\dot | \to S^t (M_s, \C)
\end{equation}
by taking the compatible family
$$
(w_n) \in \prod_{n \geq 0} A^s (\Delta^n) \otimes A^t (M_n)
$$
to the $t$-cochain
$$
\sigma \mapsto \int_{\Delta^s \times \Delta^t} w_s, \quad \sigma :
\Delta^t \to M_s.
$$

The following result is a direct consequence of Stokes' Theorem and
the compatibility conditions.

\begin{proposition}\label{A.3.3} The mapping \ref{A.3.2} is a
morphism of double complexes. \endproof
\end{proposition}

\begin{corollary}\label{A.3.4} The mapping \ref{A.3.2} induces a
morphism of complexes
$$
\int : A^\dot | M_\dot | \to S^\dot
(M_\dot,\C).
$$
\endproof
\end{corollary}

The following de~Rham theorem is a modest generalization of those of
Thom (see \cite{Sw}) and Sullivan \cite{Su}.  Proofs of statements
similar to the following appear in \cite{D1,D2}, \cite[(5.4)]{H}.

\begin{theorem}\label{A.3.5} For each simplicial manifold $M_\dot$,
the integration mapping
$$
\int : A^\dot | M_\dot | \to S^\dot
(M_\dot,
\C)
$$
is a quasi-isomorphism that induces an algebra isomorphism
$$
H^\dot (A^\dot | M_\dot |) \to
H^\dot (S^\dot (M_\dot, \C)) \cong H^\dot (|M_\dot |, \C).
$$ \endproof
\end{theorem}

\subsection{Bundles over simplicial manifolds}\label{A.4}
As is well-known, we may talk
equivalently about vector bundles or about principal bundles, and
here we choose the latter. Let $G$ be a Lie
group.  A principal $G$-bundle over the simplicial manifold
$M_\dot$ is a simplicial manifold $P_\dot$ and a morphism
$\pi : P_\dot \to M_\dot$ of simplicial manifolds satisfying:

\begin{enumerate}
\item[(i)] for each $n$, $\pi : P_n \to M_n$ is a principal
$G$-bundle;
\item[(ii)] for each morphism $f : [m] \to [n]$ of $\simp$,
$f_P$ is a morphism of $G$-bundles:
\begin{equation*}
\begin{CD}
P_n @>{f_P}>> P_m\cr
@VVV @VVV \cr
M_n @>{f_M}>> M_m\cr
\end{CD}
\end{equation*}
\end{enumerate}
The second condition guarantees that the principal
$G$-action on $P_\dot$
passes to its geometric realization:

\begin{proposition}\label{A.4.1} If $\pi : P_\dot \to M_\dot$ is a
principal $G$-bundle over a simplicial manifold, then
$| \pi| : | P_\dot | \to | M_\dot |$ is a principal $G$-bundle,
with $G$-action induced by
\begin{align*}
\Delta^n \times P_n \times G &\to \Delta^n \times P_n\cr
(t, x, g) &\mapsto(t, xg).
\end{align*}
\endproof
\end{proposition}

One can define vector bundles over simplicial manifolds similarly.

\begin{definition}\label{A.4.2} A connection $\nabla$ on the principal
$G$-bundle $P_\dot \to M_\dot$ over a simplicial manifold
$M_\dot$ is a $G$-invariant $1$-form
$$
\omega \in A^1( | P_\dot | , \Ad{\g})
$$
taking values in $\g$ (cf.~\ref{A.2.1}), the Lie algebra of
$G$, on which $G$ acts via the adjoint representation.
\end{definition}

As in the classical case, the $\g$-valued $2$-form
$$
d \omega + \frac{1}{2} [\omega, \omega] \in A^2(|P_\dot|,\Ad\g)
$$
descends to a $\g$-valued $2$-form $\Theta \in A^2(|M_\dot|, \Ad\g)$.
The $2$-form $\Theta$ is necessarily unique and is called the {\it
curvature} of $\nabla$.

We conclude this subsection with a review of the extension of
Chern-Weil theory to simplicial manifolds (cf.\ \cite{D1,D2}). We
restrict our attention  to Chern classes.

\begin{definition}\label{A.4.3} The $k{\rm th}$ {\it Chern form}
$c_k (\nabla)$ of a connection $\nabla$ on the principal
$GL_n(\C)$-bundle $P_\dot \to M_\dot$ over a simplicial manifold is:
$$
c_k (\nabla) := C_k (\Theta) \in A^{2k} | M_\dot |,
$$
where $\Theta$ is the curvature of $\nabla$ and $C_k$ is the
invariant polynomial defined in Section~\ref{1.3}.
\end{definition}

The standard arguments can be used to prove the following result (see
\cite{D2}).

\begin{theorem}\label{A.4.4} If $c_0 (\nabla), \ldots, c_n (\nabla)$
are the Chern forms of a connection on a principal $GL_n (\C)$
bundle  over a simplicial manifold, then
\begin{enumerate}
\item [(i)] each $c_k (\nabla)$ is closed;

\item [(ii)] the cohomology class of $c_k (\nabla)$ is the image
under the canonical map
$$
H^{2k} (| M_\dot |, \Z (k)) \to H^{2k} (| M_\dot|,\C)
$$
of the $k$th Chern class of $| P_\dot | \to |M_\dot |$.
\end{enumerate}
\endproof
\end{theorem}

\subsection{Cheeger-Simons classes}\label{A.5}  Suppose that $M$
is a simplicial manifold and $\Lambda$ an abelian group.  Denote
by $C^\dot (M_\dot, \Lambda)$ the cochain complex of compatible
cochains on $| M_\dot |$. This is, $C^k (M_\dot, \Lambda)$
consists of those elements $(c_n)$ of
$$
\prod_{n \geq 0} S^k (\Delta^n \times M_n)
$$
that satisfy the compatibility condition
$$
(| f | \times \id)^\ast c_n = (\id \times f_M)^\ast c_m \in
S^k (\Delta^m \times M_n)
$$
for each morphism $f : [m] \to [n]$ of $\simp$.

The following fact is standard:

\begin{proposition}\label{A.5.1} The natural inclusion
$$
C^\dot (M_\dot, \Lambda) \hookrightarrow
S^\dot (| M_\dot |, \Lambda)
$$
is a quasi-isomorphism. \endproof
\end{proposition}

A direct consequence of the compatibility condition of
\ref{A.2.1} is that integration induces a chain map
$$
I : A^\dot | M_\dot | \to C^\dot (M_\dot, \C).
$$
This is easily seen to be a quasi-isomorphism.

\begin{definition}\label{A.5.2} Suppose that $\Lambda$ is a subgroup
of $\C$.  The group of mod~$\Lambda$ {\it differential characters}
of degree $k$ of a simplicial manifold $M_\dot$ is defined by
$$
\widehat H^k (M_\dot, \Lambda) = H^{k-1} (\cone \{ A^{\geq k}
| M_\dot | \mathop {\to}\limits^{\iota_\Lambda} C^\dot
(M_\dot,\C /\Lambda) \} ),
$$
where $\iota_\Lambda$ is the composite of $I$ with the natural map
$C^\dot(M_\dot, \C) \to C^\dot (M_\dot, \C /\Lambda)$.
\end{definition}

Since the differential characters are constructed from a cone, we
have the short exact sequence
\begin{equation}\label{A.5.3}
0 \to H^{k-1} (| M_\dot |, \C /\Lambda) \to
\widehat
H^k(M_\dot, \Lambda) \to A_{cl}^k (| M_\dot |,\Lambda)
\to 0,
\end{equation}
where
$$
A_{cl}^k (| M_\dot |, \Lambda) = \{ \varphi \in A^k
|M_\dot | : d\varphi = 0 \text{ and the periods  of } \varphi
\text{ on } H_k (| M_\dot |) \text{ lie in } \Lambda \}.
$$

The construction of the Cheeger-Simons invariant $\cs_p$ given
in \ref{3.5} works equally well in the simplicial setting.
We use it to define Cheeger-Simons classes
$$
\cs_p (E_\dot, \nabla) \in
\widehat H^{2p} (M_\dot, \Z (p))
$$
for complex vector bundles $E_\dot \to M_\dot$ with connection
over simplicial manifolds.

\subsection{Simplicial varieties.}\label{5.1}
Let $Alg$ be the category whose objects are arbitrary disjoint
unions of complex algebraic manifolds and whose morphisms are
morphisms of varieties on each component.

\begin{definition}\label{5.1.1} A {\it smooth simplicial variety} is a
strict\footnote{The modifier {\it strict} means that this is the
category whose objects come {\it without} degeneracy maps.} simplicial
object of $Alg$.  The category of simplicial varieties will be denoted
by $Alg^\simp$.
\end{definition}

For later use we record the following result whose proof follows from
standard results in algebraic geometry.

\begin{proposition}\label{5.1.2}
\begin{enumerate}
\item[(a)] If $Y_\dot$ is a simplicial
variety, then there is a simplicial variety $\Ybar_\dot$
and an open immersion $Y_\dot \to \Ybar_\dot$ such that
each $\Ybar_n$ is complete and such that $\Ybar_n - Y_n$
is a normal crossing divisor in $\Ybar_n$.  Moreover, if
$f : Y_\dot \to Z_\dot$ is a morphism of simplicial
varieties, then there exist smooth completions $\Ybar_\dot$
of $Y_\dot$ and $\Zbar_\dot$ of $Z_\dot$ as above and a
morphism $\overline f: \Ybar_\dot \to \Zbar_\dot$ whose
restriction to $Y_\dot$ is $f$.
\item[(b)] If $Y^\prime_\dot$ and $Y^{\prime\prime}_\dot$
are two smooth completions as in (a) of a simplicial variety $Y_\dot$,
then there exists a third completion $\Ybar_\dot$ of $Y_\dot$
and morphisms $\Ybar_\dot \to Y^\prime_\dot$,
$\Ybar_\dot \to Y^{\prime\prime}_\dot$ such that the diagram
$$
\begin{matrix}
&&Y_\dot&&\cr
&\swarrow & \downarrow &\searrow&\cr
Y^{\prime}_\dot&\leftarrow&
\Ybar_\dot&\rightarrow&Y_\dot^{\prime\prime}
\end{matrix}
$$
commutes. \endproof
\end{enumerate}
\end{proposition}

\subsection{Complements}\label{5.2}
Because the cohomology of $B\GLn^\delta$ and other simplicial
varieties that occur naturally in the proof of Theorems~1 and~2
are not of finite type, we need to say a few words about
generalizations of Hodge theory to that case.

We will say that a simplicial manifold is of {\it finite type} if the
cohomology of its $n$-simplices is finite dimensional for each $n$.
Every simplicial
algebraic variety can be written as the direct limit of simplicial
varieties which are of finite type.  It follows that the cohomology, with
characteristic zero field coefficients, of every simplicial variety is the
inverse limit of an inverse system of mixed Hodge structures.  We view the
cohomology of a simplicial variety as a topological vector space, where
the neighbourhoods of 0 are the kernels of restriction maps to the
cohomology of simplicial varieties of finite type.  It is complete in this
topology.  We extend the notion of mixed Hodge structures to such
topological spaces.

Suppose that $\Lambda$ is a subring of $\R$ with quotient field
$\L$. Denote the category of $\Lambda$-mixed Hodge structures by
$\cH_\Lambda$. Consider the category of {\it completed $\Lambda$-mixed
Hodge structures} $\widehat\cH_\Lambda$ whose objects consist of triples
$$
(H_\Lambda, (H_\L, W_\dot), (H_\C, F^\dot)),
$$
where $H_\Lambda$ is a $\Lambda$-module; $H_\L = H_\Lambda
\otimes_\Lambda \L$ a complete topological $\L$-vectorspace with
an increasing filtration $W_\dot$ by closed subspaces; $H_\C$ the
complexification of $H_\L$ with the corresponding topology and
$F^\dot$ a decreasing filtration by closed subspaces.  In addition
we suppose  that there is a base of neighborhoods $(N_\alpha)$ of $0$
in $H_\L$ such that for each $\alpha$;
$$
\left(H_\Lambda /\left(H_\Lambda \cap N_\alpha\right),
(H_\L /N_\alpha,W_\dot),
(H_\C /N_\alpha \otimes \C, W_\dot \otimes \C,
F^\dot)\right)
$$
is a $\Lambda$-mixed Hodge structure.  A morphism
$$
(H_\Lambda, (H_\L, W_\dot), (H_\C, F^\dot)) \to
(H^\prime_\Lambda, (H^\prime_\L, W_\dot),
(H^\prime_\C, F^\dot))
$$
consists of a $\Lambda$-module homomorphism $H_\Lambda \to
H^\prime_\Lambda$ that induces continuous, filtration preserving
homomorphisms $H_\L \to H^\prime_\L$ and $H_\C \to H^\prime_\C$.

It is straightforward to show that the cohomology of every simplicial
variety inherits a functorial mixed Hodge structure, in the above sense,
from its de~Rham complex (cf.\ \cite[(8.1.19)]{De2}).

\subsection{Deligne-Beilinson (DB) cohomology}\label{5.3}
Suppose that $Y_\dot$ is a smooth simplicial variety.  By
\ref{5.1.2} we can find a smooth completion $\Ybar_\dot$
of $Y_\dot$ such that each $\Ybar_n - Y_n$ is a normal crossings
divisor $D_n$.  Denote by $I_p$ the composite
$$
F^p A^\dot (| \Ybar_\dot |\log D_\dot) \to
A^\dot |Y_\dot | \mathop {\to}\limits^I S^\dot (Y_\dot),
$$
where $F$ denotes the Hodge filtration of the de~Rham complex.

\begin{definition}\label{5.3.1} Suppose that $\Lambda$ is a subring
of $\R$ and that $p \in \N$.  The {\it Deligne-Beilinson} (or
DB) {\it cohomology} $H^\dot_\cD (Y_\dot, \Lambda (p))$ of
the simplicial variety $Y_\dot$ is the cohomology of the complex
$$
D^\dot (Y_\dot, \Ybar_\dot; \Lambda (p)) =
\cone \{ F^p A^\dot (| \Ybar_\dot | \log D_\dot )
\mathop {\to}\limits^{I_p} S^\dot (Y_\dot, \C /\Lambda
(p))\} [-1].
$$
\end{definition}

As DB-cohomology  is constructed from a cone, we have:

\begin{proposition}\label{5.3.2} For each simplicial variety
$Y_\dot$, there is natural long exact sequence
\begin{multline}
\ldots \to H^{k-1}(|Y_\dot |, \C /\Lambda (p)) \to
H_\cD^k (Y_\dot, \Lambda (p)) \cr \to F^p H^k(| Y_\dot |,\C)
\to H^k (| Y_\dot|, \C / \Lambda (p)) \to \ldots.
\end{multline}
\endproof
\end{proposition}

\begin{corollary}\label{5.3.3} For each simplicial variety
$Y_\dot$, there is a natural homomorphism
$$
H^{k-1} (|Y_\dot|, \C /\Lambda (p)) \to H_\cD^k
(Y_\dot,\Lambda (p)).
$$
\endproof
\end{corollary}

Since a point is a smooth projective variety of dimension $0$, every
simplicial set can be regarded as a simplicial variety.  Since the
Hodge filtration of the de~Rham complex of such a simplicial variety
satisfies $F^p = 0$ when $p > 0$ we have:

\begin{corollary}\label{5.3.4} If $K_\dot$ is a simplicial set and
$p > 0$, then the natural homomorphism
$$
H^{k-1} (| K_\dot |, \C / \Lambda (p)) \to H_\cD^k
(K_\dot, \Lambda (p))
$$
is an isomorphism.\endproof
\end{corollary}

\begin{corollary}\label{5.3.5} There is a natural isomorphism
$$
H^{2p}_\cD (BGL_n (\C)^\delta, \Lambda (p)) \cong
H^{2p-1} (BGL_n(\C)^\delta, \C /\Lambda (p)),
$$
where $BGL_n (\C)^\delta$ denotes the classifying space of the
general linear group endowed with the discrete topology. \endproof
\end{corollary}

\subsection{Chern classes in DB cohomology}\label{5.4}
A vector bundle of rank $r$ over a simplicial variety $Y_\dot$
is a morphism of simplicial varieties $E_\dot \to Y_\dot$ where,
for each $n$, $E_n \to Y_n$ is an algebraic
vector bundle of rank $r$. Beilinson \cite{Be} and Gillet \cite{G}
have defined Chern classes
$$
c_p^B (E_\dot) \in H_\cD^{2p} (Y_\dot, \Z (p))
$$
for vector bundles $E_\dot \to Y_\dot$ over a simplicial
varieties: Since the splitting principle holds for DB-cohomology,
the existence  of Chern classes in Deligne cohomology reduces to the
existence of first Chern classes
$$
c_1 (L_\dot) \in H^2_\cD (Y_\dot, \Z (1))
$$
for line bundles $L_\dot \to Y_\dot$.  These are constructed in
\cite[(1.7)]{Be}, in the manner of \cite{Gd}.

By a connection on a vector bundle $E_\dot \to Y_\dot$ we mean
a connection on the associated $GL_n (\C)$ bundle $P (E_\dot)
\to Y_\dot$ over the underlying simplicial manifold.

\begin{proposition}\label{5.4.1} Suppose that $Y_\dot$ is a smooth
simplicial variety and that
$$
H^{2p-1} (| Y_\dot |, \C / \Z (p)) = 0.
$$
If $\nabla$ is a connection on the vector bundle $E_\dot \to
Y_\dot$ whose $p$th Chern form satisfies
$$
c_p (\nabla) \in
F^p A^{2p} (| \Ybar_\dot | \log D_\dot),
$$
for some smooth completion $\Ybar_\dot$ of $Y_\dot$ where
$\Ybar_\dot - Y_\dot$ is a divisor with normal crossings,
then
$$
(c_p (\nabla), - \eta) \in
D^{2p} (Y_\dot, \Ybar_\dot; \Z(p))
$$
represents $c_p^B (E_\dot)$, where $\eta$ is any element of
$S^{2p-1}(Y_\dot, \C /\Z (p))$ such that
$\delta \eta = c_p (\nabla)$.
\end{proposition}

\begin{proof} The result follows from \ref{5.3.2} and the
vanishing of $H^{2p-1} (| Y_\dot|, \C /\Z (p))$ which imply that
$$
H_\cD^{2p} (Y_\dot, \Z (p)) \to F^p H^{2p} (| Y_\dot |,\C)
$$
is injective.
\end{proof}

\subsection{Multisimplicial manifolds and varieties}\label{multi}
In some constructions that one would like to do, such as
\ref{6.1.5}, one needs to work with {\it multisimplicial
varieties}.  Define the category of $0$-simplicial manifolds to be
$\cM$, the category of manifolds.  For each $n \geq 1$, we define
inductively an $n$-{\it simplicial manifold} to be a strict simplicial
object in the category of $(n-1)$-simplicial manifolds; morphisms
between two $n$-simplicial manifolds are natural transformations.
The category of $n$-simplicial manifolds and morphisms will be denoted
by $\cM^{\simp^n}$. Just as every manifold can be regarded as a
simplicial manifold, every $n$-simplicial manifold can be regarded
as an $(n+1)$ simplicial manifold by placing it in degree 0 of the
$(n+1)$-simplicial manifold, and the empty set in positive degrees.
Objects of one of the categories $\cM^{\simp^n}$
are called a {\it multisimplicial manifolds}.

The geometric realization $| M_\dot |$ of an
$n$-simplicial space
$M_\dot$ is defined inductively to be the geometric realization of
the simplicial topological space $|M|_\dot$ whose space of
$k$-simplices is the geometric realization $| M_k|$ of the
$(n-1)$-simplicial manifold $M_k$.

One can extend the definition of de~Rham complex, the de~Rham
theorem, Chern-Weil theory and Cheeger-Simons Chern classes to
multisimplicial manifolds, and the definition of Deligne-Beilinson
cohomology to multisimplicial varieties, all in a straightforward way.

\section{Additional techniques}
\label{tech}

\subsection{Killing the unipotent directions}\label{8.3}
Let $B$ be a Borel 
subgroup of $G = GL_n(\C)$. Then $G$ is a principal $B$-bundle over
the Grassmannian
$G/B$.  If one decomposes $B$ as $T \times U$, where $T$ is the maximal 
torus and
$U$ is the unipotent radical, one sees that the trouble in \ref{6.1.7}
is that
the Maurer-Cartan forms have second-order poles in the direction of $U$.
But $U$
is contractible.  One might try to replace $G$ by $G/U$, an algebraic
variety of
the same homotopy type.  We didn't see that this was helpful.

\subsection{The $Q$-filtration}\label{8.2}
Because of \ref{6.1.7}, we were led to introduce another filtration
on the de Rham complex that is coarser than the Hodge filtration,
yet induces the same filtration on cohomology.

We begin by letting $X$ be the disc $\Delta$, and $D = \{0\}$. The
logarithmic (holomorphic) de Rham complex
$$
\Omega^\dot_X(\log D) = \{\O_X \to \Omega^1_X(\log D)\}
$$
comes with its usual Hodge filtration $F$.  Let
$\Omega^\dot_X(*\! D)$ denote the larger complex of forms that
are meromorphic along $D$. We define a decreasing filtration $Q$ on
$\Omega^\dot_X(\ast D)$ by:
\begin{equation}\label{6.3.1}
Q^0 \Omega^\dot_X(\ast D) = \Omega^\dot_X(\ast D),\quad
Q^1 \Omega^\dot_X(\ast D)
= F^1\Omega^\dot_X(\log D) = \Omega^1_X(\log D)[-1],
\end{equation}
and $Q^2 \Omega^\dot_X(\ast D) = 0$. In other words, we ``discount"
a 1-form that has worse than logarithmic singularities. It is standard
that the inclusion
$\Omega^\dot_X(\log D) \subset \Omega^\dot_X(\ast D)$
is a quasi-isomorphism, and it respects the filtrations.

\begin{proposition}\label{6.3.2} The inclusion of filtered complexes
$$
(\Omega^\dot_X(\log D),F) \hookrightarrow (\Omega^\dot_X(\ast D),Q)
$$
is a filtered quasi-isomorphism.
\end{proposition}

\begin{proof} The real content of the above assertion is that
$$
Gr^0_F\Omega^\dot_X(\log D) \to Gr^0_Q\Omega^\dot_X(\ast D)
$$
is a quasi-isomorphism.  This is equivalent to the exactness of
$$
0 \to \O_X \to \O_X(\ast D)
\overset d \longrightarrow\Omega^1_X (\ast D)/\Omega^1_X(\log D)\to 0,
$$
which is elementary.
\end{proof}

Next, we consider products of the preceding.  For $X = \Delta^k$,
and $D$ the union of the coordinate axes, so that
$X - D = (\Delta^\ast)^k$, we use the tensor product of the
$Q$-filtrations of the factors to define the $Q$-filtration on
$\Omega^\dot_X(\ast D)$; and for the product of the this pair
with a polydisc $\Delta^l$, we take the tensor product of the preceding
with the Hodge filtration $F$ of $\Delta^l$.  It is an easy exercise to
check:

\begin{lemma}\label{6.3.3} Let $X$ and $D$ be as above, so
$X-D \cong (\Delta^\ast)^k \times \Delta^l$. Then the filtration
$Q$ is independent of the choice of coordinates. The inclusion
$$
(\Omega^\dot_X(\log D),F) \hookrightarrow (\Omega^\dot_X(\ast D),Q)
$$
is a filtered quasi-isomorphism.
\end{lemma}

This yields at once the following globalization:

\begin{proposition}\label{6.3.4}
Let $X$ be a complex manifold, and $D$ a divisor
with normal crossings on $X$.  Then there is a uniquely-defined
filtration $Q$ of $\Omega^\dot_X(\ast D)$, given by the above in
any system of local coordinates adapted to $D$.  Moreover, the
inclusion
$$
(\Omega^\dot_X(\log D),F) \hookrightarrow (\Omega^\dot_X(\ast D),Q)
$$
is a filtered quasi-isomorphism.
\end{proposition}

\begin{remark}\label{6.3.6}
(i) In contrast to $F$, the filtration $Q$ is {\em not}
functorial: e.g., in two variables $dz/w^2 \in Q^1$, whereas its
restriction to the diagonal is $dz/z^2 \notin Q^1$.
However, if $X$ and $D$ are as in the proposition, and $f:X'\to X$
is a morphism of complex manifolds, such that $f^\ast(D)$ is a
{\em reduced} divisor with normal crossings, then pulling back
by $f$ respects the $Q$-filtrations.

(ii) $Q$ is not multiplicative: both $dw$ and $dz/w^2$ are in $Q^1$, but
$dw \wedge dz/w^2 \notin Q^2$. However, it does induce a functorial,
multiplicative filtration on $H^\dot (X - D,\C)$,
via the quasi-isomorphism of Proposition \ref{6.3.4}.
\end{remark}

Finally, let $X_\dot$ be a smooth simplicial variety with smooth
completion $\Xbar_\dot$, such that $D_\dot = \Xbar_\dot - X_\dot$
has normal crossings, and where the face maps
preserve the $Q$ filtration. One defines in the obvious way a
subcomplex $A^\dot(|\Xbar_\dot| \ast\! D_\dot)$ of the simplicial
de~Rham complex $A^\dot(|X_\dot|)$, that contains
$A^\dot(|\Xbar_\dot|\,\log D_\dot)$.  One places the
$Q$-filtration on $A^\dot(|\Xbar_\dot|\ast\! D_\dot)$
by taking the $Q$-filtration on each
$A^\dot(|\Xbar_n|\ast \! D_n)$.
Since the Hodge filtration $F$ of
$A^\dot(|\Xbar_\dot|\,\log D_\dot)$ is defined similarly,
it follows from \ref{6.3.4} that

\begin{proposition}\label{6.3.5} Let $\Xbar_\dot$ and $D_\dot$
be as above.  Then the inclusion
$$
(A^\dot(|\Xbar_\dot|\,\log D_\dot),F)
\hookrightarrow (A^\dot(|\Xbar_\dot|\ast\! D_\dot,Q)
$$
is a filtered quasi-isomorphism.  
\end{proposition}

One can determine, at least, that the Maurer-Cartan connection forms, and then 
the Chern 
forms for the connection constructed in \ref{6.1.5}, are in $Q^1$ along the 
zero and
polar loci of the 
determinant function (recall \ref{6.1.7}).  
However, for this fact to be useful, one needs a serviceable compactification
of the simplicial 
model of $BG$.


\end{document}